\documentclass[11pt,a4paper]{amsart}
\usepackage[margin=28mm]{geometry}
\usepackage[T1]{fontenc}
\usepackage{lmodern}
\usepackage{amsmath,amssymb,mathtools,bm}
\usepackage{enumitem}
\usepackage[expansion=false]{microtype}
\usepackage{xurl}
\usepackage[unicode,pdfencoding=auto]{hyperref}
\hypersetup{colorlinks=true,linkcolor=black,citecolor=black,urlcolor=black,
  pdfauthor={Yufan Cen},
  pdftitle={Boundary remainders and optimal convergence rates for planar spectral partitions}}
\numberwithin{equation}{section}
\newtheorem{theorem}{Theorem}[section]
\newtheorem{lemma}[theorem]{Lemma}
\newtheorem{proposition}[theorem]{Proposition}
\newtheorem{corollary}[theorem]{Corollary}
\theoremstyle{definition}
\newtheorem{definition}[theorem]{Definition}
\theoremstyle{remark}

\allowdisplaybreaks[1]
\title[Boundary remainders for spectral partitions]{Boundary remainders and optimal convergence rates for planar spectral partitions}
\author{Yufan Cen}
\date{}
\begin{document}
\begin{abstract}
Let $S_m(\Omega)$ denote the infimum of the sum of the first Dirichlet eigenvalues over partitions of a bounded connected planar domain $\Omega$ into exactly $m$ nonempty connected open cells. We study the convergence of $|\Omega|S_m(\Omega)/m^2$ to its universal limit $c_*$. If the full topological boundary has zero area, we bound the absolute error by the area of its tubular neighborhood at scale $m^{-1/2}$, together with an $O(m^{-1/2})$ term. This yields rates $O_\Omega(m^{-1/2})$ on Lipschitz domains and $O_\Omega(m^{-\alpha/2})$ on domains with finitely many $C^{0,\alpha}$ boundary charts, $0<\alpha<1$, for all sufficiently large integers. The square and a fixed smooth domain have positive matching lower bounds of order $m^{-1/2}$. For each $\alpha\in(0,1)$, a fixed H\"older graph domain attains the order $m^{-\alpha/2}$ along a subsequence. Little H\"older domains satisfy the corresponding little-$o$ estimate, but no prescribed additional factor tending to zero works for the whole class, even with domain-dependent constants. An auxiliary energy with variable cell count gives the boundary estimate; a positive correction on the square and sparse graph constructions yield the lower bounds.
\end{abstract}
\maketitle

\section{Introduction and main results}\label{sec:intro}

Let $\Omega\subset\mathbb R^2$ be a bounded connected open set. For an integer $m\ge1$, consider the optimal partition energy
\begin{equation}\label{eq:partition-energy}
 S_m(\Omega)=\inf\Bigl\{\sum_{i=1}^m\lambda_1(D_i):
 (D_1,\ldots,D_m)\in\mathcal P_m(\Omega)\Bigr\}.
\end{equation}
Here $\lambda_1(D)$ is the first eigenvalue of the Dirichlet Laplacian defined on $H_0^1(D)$, and $\mathcal P_m(\Omega)$ consists of exactly $m$ nonempty connected, pairwise disjoint open sets $D_i\subset\Omega$ whose union covers $\Omega$ up to a set of Lebesgue measure zero. There are no constraints on cell areas, convexity, or interface regularity. Throughout, $|E|$ denotes two-dimensional Lebesgue measure.

Let $b=\pi j_{0,1}^2$ be the first eigenvalue of the disk of unit area, where $j_{0,1}$ is the first positive zero of the Bessel function of order zero. The Faber--Krahn inequality \cite[Theorem~3.2.1]{Henrot2006} and the Cauchy--Schwarz inequality give
\begin{equation}\label{eq:basic-growth}
 \sum_{i=1}^m\lambda_1(D_i)
 \ge b\sum_{i=1}^m|D_i|^{-1}
 \ge \frac{bm^2}{|\Omega|}.
\end{equation}
The scaling law $\lambda_1(rD)=r^{-2}\lambda_1(D)$ makes $|\Omega|S_m(\Omega)/m^2$ invariant under dilation. Let $Q=(0,1)^2$, and define
\begin{equation}\label{eq:leading-constant}
 c_*:=\inf_{n\ge1}\frac{S_n(Q)}{n^2},
 \qquad
 \delta_m(\Omega):=\frac{|\Omega|S_m(\Omega)}{m^2}-c_*.
\end{equation}
Section~\ref{upper:sec:model} proves that $c_*=\lim_{n\to\infty}S_n(Q)/n^2\in(0,\infty)$.

Geng and Lin \cite[Theorem~1.1]{GengLin2022} proved that the normalized $\ell^1$ partition energy has a positive limit independent of the shape of a smooth bounded domain. Their proof uses square replication and localization, and they note in Remark~1.1 that smoothness is not essential. We seek an estimate for $\delta_m(\Omega)$ in terms of the boundary geometry at the cell scale $m^{-1/2}$, valid for every sufficiently large integer $m$, together with fixed domains realizing the optimal orders.

Bucur, Buttazzo, and Henrot \cite{BucurButtazzoHenrot1998} proved existence results for optimal partitions in the quasi-open setting. Caffarelli and Lin \cite{CaffarelliLin2007} formulated the partition problem in terms of segregated, normalized eigenfunctions and studied the regularity of the resulting interfaces; their analysis of singularly perturbed systems \cite{CaffarelliLin2008} develops the regularity theory for segregated limits with a fixed finite number of phases. For the class in \eqref{eq:partition-energy}, our comparison arguments require connected open cells and preservation of the prescribed count. We therefore prove a spectral completion lemma and realize the auxiliary minimum by such cells. The lower-bound argument also requires compactness when the number of phases varies.

For the spectrum of a single domain, Netrusov and Safarov \cite[Theorems~1.3, 1.8, and~1.10]{NetrusovSafarov2005} obtained Dirichlet and Neumann counting estimates governed by rough boundary geometry and constructed Neumann examples attaining the H\"older order and obstructing a prescribed improvement in the little H\"older class. These examples provide a precedent for distinguishing an optimal power over a boundary class from the little-$o$ behavior on each fixed domain. In the partition problem, the construction must control the normalized Dirichlet first eigenfunctions on the optimized cells; the mass on the entire support of a cutoff determines its spectral cost.

For Lipschitz domains, Frank and Larson proved two-term asymptotics for Dirichlet Riesz means of order one \cite[Theorem~1.1]{FrankLarson2020} and, later, for Dirichlet and Neumann Riesz means of every positive order \cite{FrankLarson2025}. Their fixed-domain examples for heat traces and Dirichlet Riesz means \cite[Theorems~1.1--1.2]{FrankLarsonError2020} exclude any prescribed quantitative improvement of the remainder over the Lipschitz class. That remainder is measured after the boundary term in a spectral expansion. Here $\delta_m$ measures the error from the leading partition limit, and localization also changes the cell count. Recovering that count and retaining connected cells are additional parts of the argument.

\subsection{Boundary remainders and optimal exponents}

For $t>0$, write
\begin{equation}\label{eq:boundary-tube}
 V_\Omega(t):=\bigl|\{x\in\mathbb R^2:\operatorname{dist}(x,\partial\Omega)<t\}\bigr|.
\end{equation}
The tube includes both sides of the full topological boundary. The scale $m^{-1/2}$ is the natural cell length when the total area is fixed.

\begin{theorem}[Boundary-tube remainder]\label{thm:main-tube}
Let $\Omega\subset\mathbb R^2$ be a fixed bounded connected open set with $|\partial\Omega|=0$. There exist $C_\Omega,K_\Omega>0$ and $m_\Omega\in\mathbb N$ such that, for every integer $m\ge m_\Omega$,
\begin{equation}\label{eq:main-tube-rate}
 |\delta_m(\Omega)|\le C_\Omega
 \left[V_\Omega\!\left(K_\Omega m^{-1/2}\right)+m^{-1/2}\right].
\end{equation}
\end{theorem}

The theorem gives the universal limit whenever $|\partial\Omega|=0$. For graph boundaries, estimating the tube area gives the following rates; the graph covering and one-sided conventions are specified in Definition~\ref{def:graph-windows}.

\begin{corollary}[Rates for graph boundaries]\label{cor:regularity-rates}
Every fixed bounded connected Lipschitz domain, and hence every bounded connected $C^2$ or $C^\infty$ domain, satisfies
\[
 |\delta_m(\Omega)|=O_\Omega(m^{-1/2}).
\]
If the full boundary admits finitely many $C^{0,\alpha}$ graph charts, where $0<\alpha<1$, then
\[
 |\delta_m(\Omega)|=O_\Omega(m^{-\alpha/2}).
\]
Both estimates hold for all sufficiently large integers, with constants allowed to depend on the fixed domain.
\end{corollary}

\begin{theorem}[Optimality of the half-order rate]\label{thm:halfsharp}
There exist $\kappa_Q,C_Q>0$ such that, for all sufficiently large integers $m$,
\begin{equation}\label{eq:halfsharp-square}
 \kappa_Qm^{-1/2}\le\delta_m(Q)\le C_Qm^{-1/2}.
\end{equation}
There also exist a fixed bounded connected $C^\infty$ domain $\Omega_\infty$ and constants $\kappa_\infty,C_\infty>0$ such that, for all sufficiently large integers $m$,
\begin{equation}\label{eq:halfsharp-smooth}
 \kappa_\infty m^{-1/2}\le\delta_m(\Omega_\infty)
 \le C_\infty m^{-1/2}.
\end{equation}
\end{theorem}

\begin{theorem}[A fixed H\"older example]\label{thm:holder-sharp}
For each $0<\alpha<1$, there is a fixed bounded Jordan domain $\Omega_\alpha=\operatorname{int}\overline{\Omega_\alpha}$ whose boundary admits finitely many $C^{0,\alpha}$ graph charts and such that
\begin{equation}\label{eq:holder-positive-limsup}
 0<\limsup_{m\to\infty}m^{\alpha/2}\delta_m(\Omega_\alpha)<\infty.
\end{equation}
\end{theorem}

Theorem~\ref{thm:halfsharp} gives matching bounds for every sufficiently large integer, while Theorem~\ref{thm:holder-sharp} asserts a positive finite limsup. Both establish optimal exponents over the stated classes through fixed examples. A positive matching lower bound on every fixed Lipschitz or smooth domain remains open.

We call a domain little $C^{0,\alpha}$ if it admits finitely many boundary charts whose graph functions have uniform moduli of continuity $\omega_f(t)=o(t^\alpha)$ on their defining intervals; see Definition~\ref{def:little-holder-domain}.

\begin{theorem}[Little H\"older domains]\label{thm:little-holder}
Fix $0<\alpha<1$. Every fixed bounded connected little $C^{0,\alpha}$ domain with finitely many boundary charts satisfies
\begin{equation}\label{eq:little-pointwise}
 |\delta_m(\Omega)|=o_\Omega(m^{-\alpha/2}).
\end{equation}
Conversely, for every positive sequence $\varepsilon_m\to0$, there exists a fixed Jordan domain $\Omega_\varepsilon$ in the same class, with $\Omega_\varepsilon=\operatorname{int}\overline{\Omega_\varepsilon}$, such that
\begin{equation}\label{eq:little-no-uniform}
 \limsup_{m\to\infty}
 \frac{\delta_m(\Omega_\varepsilon)}{\varepsilon_m m^{-\alpha/2}}=+\infty.
\end{equation}
\end{theorem}

Equation~\eqref{eq:little-no-uniform} rules out an estimate $|\delta_m(\Omega)|=O_\Omega(\varepsilon_m m^{-\alpha/2})$ valid for every domain in the class, even when the constant and starting index depend on the domain.

\subsection{Elements of the proofs}
For a bounded open set $U$, a finite \emph{packing} is a family of nonempty connected, pairwise disjoint open subsets of $U$, with no covering requirement. Including the empty family, define
\[
 G(U)=\inf\Bigl\{\sum_{i=1}^n\bigl(\lambda_1(D_i)-2c_*\bigr):
 n\ge0,\ (D_i)_{i=1}^n\text{ is a packing in }U\Bigr\}.
\]
Set $B(U)=G(U)+c_*|U|$. These quantities are defined also for disconnected or empty $U$. Packings in disjoint containers can be combined, and spectral completion relates a fixed positive packing count to the partition model.

Translated cutoffs and a comparison with periodically repeated holes bound $|B(R\Omega)|$ by $C_\Omega R^2e(R)$, where $e(R)=V_\Omega(K/R)+R^{-1}\to0$. To restore the prescribed count, reserve a fixed interior square and let $T$ be the remaining region. Dilation of a nearly minimizing packing in $RT$, with count discrepancy $\Delta=n-|T|R^2$, gives
\[
 \frac{\Delta^2}{R^2}\le C_\Omega R^2e(R).
\]
The square estimate for every integer bounds the cost of repairing this discrepancy by $C_\Omega(\Delta^2/R^2+R)$. Completion therefore gives the same remainder for each prescribed integer.

The square lower bound uses reflection across its edges. At the scale where cell eigenvalues are bounded, density and stress estimates give strong compactness under one relabeling, strict local label exclusion, and uniformly bounded complete-cell diameters. The finite-phase weak-reflection theorem of Tavares and Terracini \cite[Theorem~1.1]{TavaresTerracini2012} then describes the interface of the limit. If $B(RQ)=o(R)$ along a sequence, disjoint neighborhoods of a reflected edge yield local families whose absolute defects tend to zero. Every finite deletion and insertion comparison passes to a limit containing an entire straight interface. Tangential sliding with cost $O(1/\log L)$, followed by a fixed normal deformation with a strictly negative energy change, contradicts these comparisons. The resulting linear lower bound for $B(RQ)$ gives the square rate. A complement estimate transfers it to a rounded square whose rounding scale is fixed before $m\to\infty$.

For the H\"older examples, sparse peaks of widths $w_j$ and heights $h_j=\kappa w_j^\alpha$ force zero traces along adjacent peak lines. This controls the mass on the entire cutoff support and permits removal of a core of area comparable to $h_j$. The positive area contribution dominates the perimeter and cutoff losses once the prefix is fixed and $w_j$ is chosen sufficiently small, since $w_j/h_j\to0$. Geometric decay bounds the area of all future layers, so the recursion produces one domain realizing the lower bound at $m_j\asymp w_j^{-2}$. Multiplying the heights by coefficients tending to zero gives the little H\"older example and allows the scales to accommodate any prescribed improvement sequence.

Sections~\ref{upper:sec:model}--\ref{upper:sec:classes} prove completion and the boundary estimates. Sections~\ref{sec:compactness}--\ref{fullsquare:sec} establish the square lower bound and the smooth example; Sections~\ref{holder:sharp-section}--\ref{holder:little-section} contain the H\"older constructions. Appendix~\ref{app:heat} proves the heat-semigroup domination used for low-energy eigenfunctions on arbitrary open cells.

\smallskip
\noindent\textbf{Use of generative AI.} OpenAI's GPT-6-Astra, GPT-6-Sol, and GPT-6.1-Sol were used during this work.
Building on the author's spectral completion lemma (Lemma~\ref{upper:lem:completion}), discussions with AI led to the auxiliary-energy approach in Sections~\ref{upper:sec:squareupper}--\ref{upper:sec:classes}.
Through continued interaction with the author, AI developed the sharp $m^{-1/2}$ lower-bound proofs for the square and smooth-domain examples suggested by the author (Sections~\ref{sec:compactness}--\ref{fullsquare:sec}).
It also proposed and developed the H\"older and little H\"older constructions (Sections~\ref{holder:sharp-section}--\ref{holder:little-section}) and assisted with literature searches, drafting, and revision.
The author independently checked the main proofs.
Additional audits were conducted separately using DeepSeek-V4.1-Flash and GPT-6-Astra.
The author takes full responsibility for the final manuscript.

\section{The partition model and the universal constant}\label{upper:sec:model}

The partition class $\mathcal P_m(\Omega)$ and the energy $S_m(\Omega)$ are as in \eqref{eq:partition-energy}. For a bounded open set $D$,
\[
 \lambda_1(D)=\inf_{0\ne u\in H_0^1(D)}
 \frac{\int_D|\nabla u|^2}{\int_Du^2}.
\]
First eigenfunctions are taken nonnegative, normalized in $L^2$ on their complete cells, and extended by zero. We use $S_m(r\Omega)=r^{-2}S_m(\Omega)$.

\subsection{Completion without changing the cell count}
\begin{lemma}[Spectral completion]\label{upper:lem:completion}
Let $\Omega$ be a bounded connected open set, let $n\ge1$, and let $U_1,\ldots,U_n\subset\Omega$ be nonempty pairwise disjoint open sets. For every $\varepsilon>0$, there exists $(D_i)_{i=1}^n\in\mathcal P_n(\Omega)$ such that
\[
 \sum_i\lambda_1(D_i)\le\sum_i\lambda_1(U_i)+\varepsilon.
\]
\end{lemma}
\begin{proof}
If $n=1$, take $D_1=\Omega$ and use domain monotonicity. Suppose $n\ge2$. For each $i$, choose a nonzero smooth compactly supported function in $U_i$ with Rayleigh quotient less than $\lambda_1(U_i)+\varepsilon/(2n)$. An open set has at most countably many connected components, and its $H_0^1$ space is the Hilbert direct sum of the spaces on those components. Since a Rayleigh quotient is the mass-weighted average of the component quotients, there is a component $W_i$ with
\[
 \lambda_1(W_i)\le\lambda_1(U_i)+\frac{\varepsilon}{2n}.
\]
Choose a new smooth compactly supported test function in $W_i$ sufficiently close to its first eigenvalue. Finitely many small balls and polygonal paths inside $W_i$ connect this support within a connected open set $V_i\Subset W_i$, with
\[
 \lambda_1(V_i)\le\lambda_1(U_i)+\frac{\varepsilon}{n}.
\]
The sets $K_i=\overline V_i$ are compact and connected, and their pairwise distances are positive.

Let $d_i(x)=\operatorname{dist}_\Omega(x,K_i)$ be the intrinsic distance measured by lengths of curves in $\Omega$. Connected open sets are polygonally connected, so $d_i$ is finite and is $1$-Lipschitz on every ball contained in $\Omega$. Set $\delta=\min_{i\ne j}\operatorname{dist}(K_i,K_j)>0$. Choose weights $|w_i|<\delta/4$ such that every set
\[
 \{d_i-d_j=w_j-w_i\},\qquad i\ne j,
\]
has measure zero. This is possible because a measurable real function has at most countably many level sets of positive measure; it suffices to avoid countably many hyperplanes in the space of weights. Define
\[
 D_i=\{x\in\Omega:d_i(x)+w_i<d_j(x)+w_j\text{ for all }j\ne i\}.
\]
For $x\in K_i$, one has $d_i(x)=0$ and $d_j(x)\ge\delta$, hence $K_i\subset D_i$. The sets $D_i$ are open and disjoint, and their complement is contained in a finite union of the null level sets just chosen.

It remains to prove connectivity. Fix $x\in D_i$, and let $\eta>0$ be the minimum of its strict advantages over the other indices. Choose a curve from $K_i$ to $x$ of length $L<d_i(x)+\eta/2$. If $z$ is on this curve and the remaining length to $x$ is $\ell$, then
\[
 d_i(z)\le L-\ell,\qquad d_j(z)\ge d_j(x)-\ell.
\]
These inequalities show that the entire curve lies in $D_i$. Every point of $D_i$ can therefore be joined to $K_i$. Path components of an open set are relatively open and closed; the connected set $K_i$ lies in a single such component. Hence $D_i$ is connected. Finally, domain monotonicity gives $\lambda_1(D_i)\le\lambda_1(V_i)$, and summing proves the lemma.
\end{proof}

Consequently, for each fixed $n\ge1$, the infimum over complete partitions equals the infimum over packings with $n$ cells. Placing any finite number of small squares inside $\Omega$ and applying the lemma also shows that every $S_m(\Omega)$ is finite.

\subsection{The square constant}
For $c_*$ as in \eqref{eq:leading-constant}, Faber--Krahn and the one-cell partition give
\[
 b\le c_*\le2\pi^2.
\]
To prove convergence to $c_*$, we use the square replication argument of \cite[Lemma~2.1 and the first step of the proof of Theorem~1.1]{GengLin2022}, together with completion.

Fix $n$, and for $m\ge n$ put $k=\lceil\sqrt{m/n}\rceil$ and $N=nk^2\ge m$. Replicate a nearly minimizing $n$-partition in the $k^2$ subsquares. The infimum of the resulting spectral sums is $k^4S_n(Q)$. Delete $N-m$ cells and apply Lemma~\ref{upper:lem:completion} to obtain
\[
 \frac{S_m(Q)}{m^2}\le\frac{S_n(Q)}{n^2}\left(\frac Nm\right)^2.
\]
For fixed $n$, $N/m\to1$. Taking the infimum over $n$ and using the defining lower bound for $c_*$ proves
\begin{equation}\label{upper:eq:squareconstant}
 \lim_{m\to\infty}\frac{S_m(Q)}{m^2}=c_*.
\end{equation}
\subsection{Elementary packing bounds}
For the auxiliary energy $G$ defined in Section~\ref{sec:intro}, Faber--Krahn gives $-c_*^2|U|/b\le G(U)\le0$ on every bounded open set $U$, with $G(\varnothing)=0$. Indeed, a packing with $n>0$ cells has spectral sum at least $bn^2/|U|$, and minimizing this quadratic lower bound over real $n$ proves the claim. Writing $Q_L=LQ$, completion and $S_n(Q)\ge c_*n^2$ give $G(Q_L)\ge-c_*L^2$ for every $L>0$.

\section{Boundary-tube estimates and prescribed cell counts}\label{upper:sec:squareupper}

\subsection{Translated cutoffs}
\begin{lemma}[Translation and cutoff]\label{upper:lem:translation}
Let $U\subset\mathbb R^2$ be bounded and open. If $\chi\in H_0^1(U)$ and $0\le\chi\le1$, then
\begin{equation}\label{upper:eq:translation}
 G(U)\le\int_U|\nabla\chi|^2-c_*\int_U\chi^2.
\end{equation}
\end{lemma}
\begin{proof}
First suppose $\chi\in C_c^\infty(U)$. By the definition of $c_*$, choose an $n$-partition of the unit square with spectral sum at most $(c_*+\epsilon)n^2$. After dilation by $\sqrt n$, its area and cell count are both $n$, and its spectral sum is $E\le(c_*+\epsilon)n$. Repeat this template periodically, translate it by $z$, and assign independent labels to the copies. Let $u_{i,z}$ be the normalized first eigenfunction on each complete cell, extended by zero. The index $i$ includes both the template label and the copy. Only finitely many labels have nonzero cutoffs for each translation. Existence of a first eigenfunction on a bounded open cell follows from compactness of $H_0^1$ into $L^2$ and requires no boundary regularity.

For $a_i=\|\chi u_{i,z}\|_2^2>0$, let $q_i=\int|\nabla(\chi u_{i,z})|^2/a_i$. The cutoff belongs to $H_0^1(D_{i,z}\cap U)$: approximate $u_{i,z}$ by smooth compactly supported functions in $D_{i,z}$ and multiply by $\chi$. A component quotient is a mass-weighted average, so one connected component has first eigenvalue at most $q_i$, or arbitrarily close to it. Retain only labels for which $q_i<2c_*$. Since $a_i\le1$ and the retained summands are negative,
\[
 G(U)\le\sum_{q_i<2c_*}(q_i-2c_*)
 \le\sum_{q_i<2c_*}a_i(q_i-2c_*)
 \le\sum_i a_i(q_i-2c_*).
\]
Zero-mass labels are omitted. Testing the weak eigenvalue equation with $\chi^2u_{i,z}$ yields
\[
 \int|\nabla(\chi u_{i,z})|^2
 =\lambda_i\int\chi^2u_{i,z}^2+\int|\nabla\chi|^2u_{i,z}^2.
\]
Let $Z$ be a period box of area $n$. Complete-cell normalization and Tonelli's theorem give, for almost every $x$,
\[
 \frac1n\int_Z\sum_i u_{i,z}(x)^2\,dz=1,\qquad
 \frac1n\int_Z\sum_i\lambda_i u_{i,z}(x)^2\,dz=E/n.
\]
Averaging the preceding bound over $z$ therefore gives
\[
 G(U)\le\int_U|\nabla\chi|^2+(E/n-2c_*)\int_U\chi^2
 \le\int_U|\nabla\chi|^2-(c_*-\epsilon)\int_U\chi^2.
\]
Let $\epsilon\downarrow0$.

For a general $[0,1]$-valued $H_0^1$ function, first take a strong $H^1$ approximation by smooth compactly supported functions and then truncate by $s\mapsto\min(1,\max(0,s))$. Truncation is continuous in $H^1$. Convolution with a nonnegative smooth kernel of sufficiently small support gives a smooth approximation with the same range and with support compactly contained in $U$. Passing to the limit in the final inequality, which contains no eigenfunction weights, proves the result.
\end{proof}

\subsection{Recovering every integer count in a square}
A near minimizer for $G(Q_R)$ gives a count close to $R^2$; a rectangular construction corrects its discrepancy.

\begin{proposition}[Square estimate for every integer]\label{upper:prop:squareupper}
There is a finite constant $K_0$ such that, for every integer $m\ge1$,
\begin{equation}\label{upper:eq:squareupper}
 c_*m^2\le S_m(Q)\le c_*m^2+K_0m^{3/2}.
\end{equation}
\end{proposition}
\begin{proof}
Set $e_n=S_n(Q)-c_*n^2\ge0$ and $S_0=e_0=0$. Completion and scaling give the exact integer envelope
\[
 B(Q_R)=\min_{n\ge0}\frac{e_n+c_*(n-R^2)^2}{R^2}.
\]
This is minimization over integer values of $n$. The quadratic term ensures attainment and does not require attainment of any $S_n(Q)$.

For $R\ge2$, apply Lemma~\ref{upper:lem:translation} with $\chi=\min(1,\operatorname{dist}(x,Q_R^{\complement}))$. After centering,
\[
 B(Q_R)\le\int_{Q_R}\bigl[|\nabla\chi|^2+c_*(1-\chi^2)\bigr].
\]
The integrand is supported in a boundary strip of width $1$, is at most $1+c_*$, and the strip has area at most $4R$. Thus $0\le B(Q_R)\le4(c_*+1)R$. For $R<2$, the empty packing gives $B(Q_R)\le c_*R^2\le2c_*R$. Consequently, there is $K$ such that every real $x>0$ admits an integer $n\ge0$ satisfying
\begin{equation}\label{upper:eq:witness}
 e_n+c_*(n-x)^2\le Kx^{3/2}.
\end{equation}
With $a=\sqrt{K/c_*}$, this implies $|n-x|\le ax^{3/4}$.

We shall place square templates in rectangles. Take a nearly minimizing $n$-partition of the unit square and its normalized first eigenfunctions $u_i$. Write
\[
 X=\sum_{i=1}^n\int|\partial_1u_i|^2,\qquad
 Y=\sum_{i=1}^n\int|\partial_2u_i|^2,
 \qquad X+Y\le S_n(Q)+\epsilon.
\]
Map either this partition or its $90$-degree rotation to a rectangle of width $w$ and height $h$. The two spectral upper bounds are $X/w^2+Y/h^2$ and $Y/w^2+X/h^2$. Choosing the better of these two competitors gives
\begin{equation}\label{upper:eq:rectangle}
 E_{w,h}\le\frac{S_n(Q)+\epsilon}{wh}\Phi(w/h),\qquad
 \Phi(p)=\frac{p+p^{-1}}2,\quad \Phi(p)-1=\frac{(p-1)^2}{2p}.
\end{equation}
Only the two numerical bounds have been averaged; no equality between $X$ and $Y$ is assumed.

For sufficiently large $m$, choose $n$ from \eqref{upper:eq:witness} with $x=m-2am^{3/4}$. Set $d=m-n$, $t=d/m$, and $H=n/m$. Then
\[
 am^{3/4}\le d\le3am^{3/4},\quad n\ge m/2,
 \quad t=O(m^{-1/4}),\quad e_n=O(m^{3/2}).
\]
Place the $n$-cell template in the main rectangle of width $1$ and height $H$. Equation~\eqref{upper:eq:rectangle}, with the template error subsequently sent to zero, gives
\begin{align*}
 E_{\rm main}
 &\le\tfrac12(1+H^{-2})(c_*n^2+e_n)\\
 &=c_*mn+\tfrac{c_*}{2}d^2+O(m^{3/2})
 =c_*mn+O(m^{3/2}).
\end{align*}
The remaining strip has height $t$. Put $x_0=d^2/m\asymp m^{1/2}$ and apply \eqref{upper:eq:witness} again to obtain
\[
 n_0=x_0+O(x_0^{3/4}),\qquad e_{n_0}=O(x_0^{3/2}).
\]
Let $b_0=\lceil d/n_0\rceil$, $q=b_0n_0=d+r$, and $0\le r<n_0$. Divide the strip into $b_0$ rectangles, each of width $1/b_0$ and height $t$, and place an $n_0$-cell template in each. Since
\[
 p=b_0t=x_0/n_0+O(t)=1+O(x_0^{-1/4}+t),
\]
we have $\Phi(p)=1+O(x_0^{-1/2}+t^2)$. Using $\Phi(p)=\Phi(p^{-1})$, we obtain
\[
 E_{\rm strip}\le\frac{q^2}{t}
 \left(c_*+\frac{e_{n_0}}{n_0^2}\right)\Phi(p).
\]
Now $x_0=dt$ and $q=d+O(x_0)$ imply $q^2/t=md+O(d^2)$, while $e_{n_0}/n_0^2=O(x_0^{-1/2})$. Hence
\begin{align*}
 E_{\rm strip}
 &\le\bigl(md+O(d^2)\bigr)\bigl(c_*+O(x_0^{-1/2})\bigr)\\
 &\qquad\times\bigl[1+O(x_0^{-1/2}+t^2)\bigr]
 \le c_*md+O(m^{3/2}),
\end{align*}
because
\[
 mdx_0^{-1/2}=m^{3/2},\qquad
 mdt^2=d^3/m=O(m^{5/4}),
\]
and $d^2=O(m^{3/2})$.

The two regions contain $m+r$ cells. Delete $r$ cells, which cannot increase the spectral sum, and complete the remaining packing. For each fixed $m$, let the finitely many template errors tend to zero. This gives the upper bound for large $m$, and enlarging the constant covers the remaining integers. The lower bound is the definition of $c_*$.
\end{proof}

\subsection{The boundary-tube estimate on a general domain}\label{upper:sec:tube}
The scale factor remains inside $V_\Omega$; no doubling assumption is made.

\begin{proof}[Proof of Theorem~\ref{thm:main-tube}]
For a bounded open set $U$, let $I_U(t)$ be the area of the points in $U$ at distance less than $t$ from $U^{\complement}$, and let $O_U(t)$ be the area of the points outside $\overline U$ at distance less than $t$ from $\overline U$.

The function $\chi=\min(1,\operatorname{dist}(x,U^{\complement})/t)$ belongs to $H_0^1(U)$. To see this directly, extend $d(x)=\operatorname{dist}(x,U^{\complement})$ by zero to the plane and set $\chi_\epsilon=\min(1,(d-\epsilon)_+/t)$. Boundedness of $U$ makes its support compact in $U$, and
\[
 \nabla\chi_\epsilon=t^{-1}\nabla d\,\mathbf1_{\{\epsilon<d<\epsilon+t\}}.
\]
The gradient of a Sobolev function vanishes almost everywhere on each level set. Dominated convergence therefore gives strong $H^1$ convergence, and smooth approximation proves the claim. Lemma~\ref{upper:lem:translation} yields
\[
 B(U)\le(c_*+t^{-2})I_U(t).
\]

For the lower bound, assume $|\partial U|=0$. Place $U$ inside a square of side $L$ at distance greater than $2t$ from its edges. Put $N^2$ copies in an $N$ by $N$ array forming a large square $W_N$, and let
\[
 K_N=W_N\setminus\bigcup_j\overline U_j.
\]
Use $\chi_N=\min(1,\operatorname{dist}(x,K_N^{\complement})/t)$. The exterior layers of the holes are disjoint and do not meet the outer square layer. The latter has area at most $4NLt$, and each hole layer has area $O_U(t)$. Hence
\[
 \int_{K_N}\bigl[|\nabla\chi_N|^2+c_*(1-\chi_N^2)\bigr]
 \le N^2(c_*+t^{-2})O_U(t)+O_{U,t}(N).
\]
Combining finite nearly minimizing packings in the separate containers and sending their errors to zero gives $G(W_N)\le N^2G(U)+G(K_N)$. Use the square lower bound, the cutoff upper bound on $K_N$, and the identity $|K_N|=|W_N|-N^2|U|$. Divide by $N^2$ and send $N\to\infty$ to conclude
\begin{equation}\label{upper:eq:grandtube}
 -(c_*+t^{-2})O_U(t)\le B(U)\le(c_*+t^{-2})I_U(t).
\end{equation}
Every comparison uses a finite packing. The zero-area boundary assumption enters the area identity through $|\overline U|=|U|$; there is no residual area term after centering.

We record a consequence used later. If $\partial U$ is a rectifiable Jordan curve of length $P_U$, divide its arc-length interval into at most $P_U/t+1$ pieces of length at most $t$. Each piece lies in a radius-$t$ ball centered at an endpoint, so its $t$-neighborhood lies in the corresponding radius-$2t$ ball. Thus
\begin{equation}\label{upper:eq:rectifiabletube}
 V_U(t)\le4\pi(P_Ut+t^2).
\end{equation}
A rectifiable curve has zero area. Applying \eqref{upper:eq:grandtube} to $RU$ with $t=1$ gives
\begin{equation}\label{upper:eq:perimetergrand}
 |G(RU)+c_*|U|R^2|
 \le4\pi(c_*+1)(P_U R+1),\qquad R>0.
\end{equation}
The constant $C_0=4\pi(c_*+1)$ is independent of $U$.

Write $A=|\Omega|$, and fix a square $C\Subset\Omega$ of area $v>0$ such that $T=\Omega\setminus\overline C$ has area $A_T=A-v>0$. The set $T$ may be disconnected, and $V_T(t)\le V_\Omega(t)+C_Ct$ for all sufficiently small $t$. Compactness and zero area of $\partial\Omega$ imply $V_\Omega(t)\to0$. Scaling \eqref{upper:eq:grandtube} with physical tube width $1$, for $U=\Omega,T$ and sufficiently large $r$ we obtain
\begin{equation}\label{upper:eq:grandmodulus}
 |G(rU)+c_*|U|r^2|\le C_\Omega r^2[V_\Omega(1/r)+1/r].
\end{equation}
Choose a finite packing in $RT$ with auxiliary energy at most $G(RT)+1$, and delete every cell with first eigenvalue at least $2c_*$. Denote its remaining count and spectral sum by $n$ and $E$. Then
\[
 G(RT)\le E-2c_*n\le G(RT)+1.
\]
Faber--Krahn gives $n\le2c_*A_TR^2/b$, while \eqref{upper:eq:grandmodulus} gives $E-2c_*n\le-\tfrac12c_*A_TR^2$ for large $R$. Thus $n\ge A_TR^2/4$. In particular,
\[
 s=\sqrt{n/A_T},\qquad R/2\le s\le\sqrt{2c_*/b}\,R.
\]
Fix $K\ge2$ and put $e(R)=V_\Omega(K/R)+1/R\to0$. Monotonicity of $V_\Omega$ bounds the errors at both scales $R,s$ by $C_\Omega R^2e(R)$. Dilating the packing from $RT$ to $sT$ yields
\[
 R^2E/s^2-2c_*n\ge G(sT)\ge-c_*n-C_\Omega R^2e(R).
\]
Multiply by $s^2/R^2$ and combine with the upper bound at scale $R$ to get
\[
 E\ge\frac{c_*n^2}{A_TR^2}-C_\Omega R^2e(R),\qquad
 E-2c_*n\le-c_*A_TR^2+C_\Omega R^2e(R).
\]
Completing the square, with $\Delta=n-A_TR^2$, gives
\begin{equation}\label{upper:eq:countsquare}
 \frac{c_*\Delta^2}{A_TR^2}\le C_\Omega R^2e(R),\qquad
 \Delta=o(R^2).
\end{equation}
Indeed, subtraction produces $c_*[n^2/(A_TR^2)-2n+A_TR^2]\le C_\Omega R^2e(R)$, so $|\Delta|/R^2\le C_\Omega\sqrt{e(R)}\to0$.

Given an integer $m$, set $R=\sqrt{m/A}$ and $k=m-n=vR^2-\Delta$. For large $m$, $k$ is a positive integer and $k/R^2\to v$. Given $\eta>0$, Proposition~\ref{upper:prop:squareupper} supplies a complete $k$-partition of $RC$ with spectral sum $E_C$ such that
\begin{align*}
 E_C-2c_*k
 &\le-c_*vR^2+\frac{c_*\Delta^2}{vR^2}
 +\frac{K_0k^{3/2}}{vR^2}+\eta\\
 &\le-c_*vR^2+C_\Omega R^2e(R)+\eta.
\end{align*}
The term involving $k^{3/2}$ is $O_\Omega(R)$ and is absorbed because $e(R)\ge1/R$. Combine this family with the packing in $RT$ to obtain exactly $m$ cells. Completion in the connected set $R\Omega$ costs at most another $\eta$. Sending $\eta\downarrow0$ for fixed $m$ yields
\[
 S_m(R\Omega)\le c_*m+C_\Omega R^2e(R).
\]
Conversely, every complete $m$-partition is admissible for $G(R\Omega)$, so \eqref{upper:eq:grandmodulus} gives
\[
 S_m(R\Omega)\ge c_*m-C_\Omega R^2e(R).
\]
Divide by $m=AR^2$, use $S_m(R\Omega)/m=AS_m(\Omega)/m^2$, and absorb the fixed scale $\sqrt A$ into the constants. This proves \eqref{eq:main-tube-rate}.
\end{proof}

\section{Consequences for Lipschitz and H\"older boundaries}\label{upper:sec:classes}

\begin{definition}[Boundary regularity]\label{def:graph-windows}\label{def:little-holder-domain}
Let $\Omega\subset\mathbb R^2$ be bounded, open, and connected. A finite continuous graph covering of $\partial\Omega$ consists of rigid coordinate changes $\Phi_j$, bounded open intervals $I_j,J_j$, and functions $f_j\in C(I_j)$ with $f_j(I_j)\subset J_j$, such that
\[
 \begin{gathered}
 \partial\Omega\subset\bigcup_{j=1}^N W_j,\qquad
 W_j=\Phi_j^{-1}(I_j\times J_j),\\
 \Phi_j(\partial\Omega\cap W_j)=\{(s,f_j(s)):s\in I_j\}.
 \end{gathered}
\]
The covering is $C^{0,\alpha}$, $0<\alpha\le1$, if each $f_j\in C^{0,\alpha}(I_j)$. For $0<\alpha<1$, it is \emph{little} $C^{0,\alpha}$ if
\[
 \omega_{f_j,I_j}(r)=o(r^\alpha)\quad(r\downarrow0)
 \quad\text{for every }j,
\]
where $\omega_{f,I}(r)=\sup\{|f(s)-f(t)|:s,t\in I,\ |s-t|\le r\}$. The seminorms and moduli are taken over the entire defining intervals.

A covering is \emph{one-sided} if, after possible vertical reflections in the charts,
\[
 \Phi_j(\Omega\cap W_j)
 =\{(s,t)\in I_j\times J_j:t>f_j(s)\}.
\]
We call $\Omega$ a (little) $C^{0,\alpha}$ graph domain if it admits a corresponding one-sided covering. A Lipschitz domain means a domain with a one-sided $C^{0,1}$ covering.
\end{definition}

A finite graph covering implies $|\partial\Omega|=0$ by Fubini's theorem. The upper bounds require only this covering; the examples in Sections~\ref{holder:sharp-section} and~\ref{holder:little-section} also satisfy the one-sided condition.

\begin{proof}[Proof of Corollary~\ref{cor:regularity-rates}]
By compactness of $\partial\Omega$, choose open sets $\widetilde W_j$ with $\overline{\widetilde W_j}\Subset W_j$ that still cover $\partial\Omega$, and choose
\[
 0<t_0<\min_j\operatorname{dist}
       (\overline{\widetilde W_j},W_j^{\complement}).
\]
For $t<t_0$, each tube point lies in the same larger chart as a boundary point at distance less than $t$ in a smaller chart. In chart $j$, its vertical distance from the graph at its own horizontal coordinate is at most $t+\omega_{f_j,I_j}(t)$. Hence
\begin{equation}\label{upper:eq:charttube}
 V_\Omega(t)\le2\sum_{j=1}^N |I_j|
        \bigl[t+\omega_{f_j,I_j}(t)\bigr],
 \qquad 0<t<t_0.
\end{equation}
For Lipschitz charts this is $O_\Omega(t)$; for $C^{0,\alpha}$ charts it is $O_\Omega(t^\alpha)$. Theorem~\ref{thm:main-tube} gives the stated rates. The smooth cases follow from the Lipschitz case.
\end{proof}

If $\omega_{f_j,I_j}\le\omega$ for every $j$, where $\omega$ is nondecreasing and $\omega(t)\to0$ as $t\downarrow0$, \eqref{upper:eq:charttube} gives $V_\Omega(t)\le C_\Omega[t+\omega(t)]$, hence
\begin{equation}\label{upper:eq:graphmodulus}
 |\delta_m(\Omega)|\le C_\Omega
 [m^{-1/2}+\omega(K_\Omega m^{-1/2})]
\end{equation}
for all sufficiently large integers, with constants and starting index depending on the fixed domain. No concavity or doubling assumption on $\omega$ is used; the factor $K_\Omega$ remains inside the modulus.

\section{Uniform estimates and compactness for minimizing configurations}\label{sec:compactness}

Throughout this section set $\Lambda=2c_*$.

\subsection{Exact minimizers and estimates independent of the number of phases}
We use the scalar density and stress to obtain a gradient bound for finite or countable segregated families.

\begin{lemma}[A gradient bound from density and stress]\label{fullsquare:gradient}
On $B_{2r}$, let $U=(u_i)$ be a finite or countable nonnegative segregated family in $H^1_{\rm loc}(B_{2r};\ell^2)$. Suppose $0\le\mu_i\le\Lambda$ and $\rho=\sum_i u_i^2\le M^2$. Assume that
\begin{equation}\label{fullsquare:mult}
 \int\nabla u_i\cdot\nabla(\phi u_i)=\mu_i\int\phi u_i^2
\end{equation}
for every smooth compactly supported $\phi$, and that
\[
 \operatorname{div}T=0,\qquad
 T=(g-P)I-2\sum_i\nabla u_i\otimes\nabla u_i,
 \quad g=\sum_i|\nabla u_i|^2,\quad P=\sum_i\mu_i u_i^2.
\]
Then
\[
 \operatorname*{ess\,sup}_{B_r}g\le C(\Lambda r^2)M^2/r^2,
\]
with a constant independent of the number of labels.
\end{lemma}
\begin{proof}
The case $M=0$ is immediate. Otherwise rescale to $r=M=1$, and denote the rescaled spectral bound by $K$. Summing \eqref{fullsquare:mult} with $\phi=\eta^2$ gives $\int_{B_{7/4}}g\le C(1+K)$. Set
\[
 \partial=(\partial_x-\mathrm i\partial_y)/2,\qquad
 \bar\partial=(\partial_x+\mathrm i\partial_y)/2,\qquad
 \mathcal Q=\sum_i(u_{i,x}-\mathrm i u_{i,y})^2.
\]
Segregation and the almost-everywhere vanishing of a Sobolev gradient on the zero set imply $|\mathcal Q|=g$. The stress equation becomes
\[
 \bar\partial\mathcal Q=-\partial P,
 \qquad |\nabla P|\le2K\sqrt\rho\sqrt g\le2K\sqrt{|\mathcal Q|}.
\]
In particular,
\[
 \int_{B_{7/4}}|\partial P|^2
 \le\int_{B_{7/4}}|\nabla P|^2
 \le4K^2\int_{B_{7/4}}g
 \le CK^2(1+K).
\]
Thus $\mathcal Q\in L^1$ and $\partial P\in L^2$, with bounds depending only on $K$.

Choose a cutoff of the right-hand side equal to one on a smaller ball and take its Cauchy transform. Its first derivatives are bounded in $L^2$, and its zeroth-order term is locally in $L^2$, since the Cauchy kernel is locally integrable. The difference from $\mathcal Q$ is distributionally holomorphic on the smaller ball, hence holomorphic; its $L^1$ interior estimates control all interior derivatives. We obtain
\[
 \|\mathcal Q\|_{W^{1,2}(B_{4/3})}\le C_K.
\]
The planar Sobolev embedding gives $\mathcal Q\in L^4$, and the displayed inequality then gives $\nabla P\in L^8$. Repeating the interior Cauchy--Riemann estimate in $L^4$ yields
\[
 \|\mathcal Q\|_{W^{1,4}(B_{5/4})}\le C_K.
\]
Morrey's embedding completes the argument.

The singular-integral input is the Cauchy-transform derivative with multiplier
\[
 m(\xi)=\frac{\xi_1-\mathrm i\xi_2}{\xi_1+\mathrm i\xi_2},\qquad \xi\ne0.
\]
It is smooth away from zero and homogeneous of degree zero. The multiplier theorem in \cite[Chapter~VI, \S7.6, p.~263]{Stein1993} gives the two $L^p$ estimates used above, at $p=2,4$.
\end{proof}

For open $O\subset\mathbb R^2$, whenever $u\in H^1_{\rm loc}(O)$ admits a continuous representative $\widetilde u$, we use that unique representative and define its canonical positivity set by
\[
 D_O(u):=\{x\in O:\widetilde u(x)>0\}.
\]
\begin{lemma}[Zero trace on positivity components]\label{fullsquare:positivity-trace}
Let $O\subset\mathbb R^2$ be open, and let $u\in H_0^1(O)$ be nonnegative with a continuous representative. If $D=D_O(u)$, then the zero extension of $u$ belongs to $H_0^1(D)$. For each connected component $C$ of $D$, the zero extension of $u\mathbf1_C$ belongs to $H_0^1(C)$.
\end{lemma}
\begin{proof}
Choose nonnegative $v_n\in C_c^\infty(O)$ converging strongly to $u$ in $H^1$. Such an approximation follows by compactly supported approximation, taking positive parts, and local smoothing preserving nonnegativity. Continuity of the positive-part map in $H^1$ gives
\[
 w_n:=\min(u,v_n)=u-(u-v_n)_+\longrightarrow u
 \quad\text{strongly in }H^1.
\]
For fixed $n$ and $\epsilon>0$, put $w_{n,\epsilon}=(w_n-\epsilon)_+$. Its support is contained in
\[
 \operatorname{supp}v_n\cap\{u\ge\epsilon\}\Subset D.
\]
Smoothing inside $D$ shows that $w_{n,\epsilon}\in H_0^1(D)$. First let $\epsilon\downarrow0$ and then $n\to\infty$ to prove the first assertion.

Choose a sequence $z_k$ among these truncations such that $z_k\to u$ strongly in $H^1$ and each support is compact in $D$. Such a support meets only finitely many components: cover it by balls compactly contained in $D$ and take a finite subcover. Its intersection with each component is compact in that component. A smooth cutoff therefore gives $z_k\mathbf1_C\in H_0^1(C)$ with weak gradient $\mathbf1_C\nabla z_k$. Applying the same argument to $z_k-z_l$ yields
\[
 \|z_k\mathbf1_C-z_l\mathbf1_C\|_{H^1(\mathbb R^2)}
 \le\|z_k-z_l\|_{H^1(\mathbb R^2)}.
\]
The restricted sequence is Cauchy in $H_0^1(C)$, and its $L^2$ limit is $u\mathbf1_C$. This identifies the required weak gradient without assuming that $\partial C$ has zero measure.
\end{proof}

\begin{lemma}[Exact realization by connected open cells]\label{fullsquare:existence}
For every bounded open set $O$, $G(O)$ is attained by a finite packing of connected open cells. An attaining family can be chosen whose cells are the canonical positivity sets of nonnegative segregated normalized first eigenfunctions and whose eigenvalues satisfy $\lambda_i<\Lambda$. If $O$ is connected and $G(O)<0$, these positivity sets cover $O$ almost everywhere.
\end{lemma}
\begin{proof}
For $n\ge1$, set
\[
 \mathcal A_n(O)=\left\{(u_i)_{i=1}^n\in H_0^1(O)^n:
 u_i\ge0,\ \|u_i\|_2=1,\ u_i u_k=0\ \text{a.e. for }i\ne k\right\},
\]
and include $\mathcal A_0(O)=\{\varnothing\}$. Write
\[
 \mathcal E(U)=\sum_{i=1}^n\left(\int_O|\nabla u_i|^2-\Lambda\right),
 \qquad G_f=\inf_{n\ge0}\inf_{U\in\mathcal A_n(O)}\mathcal E(U).
\]
Normalized first eigenfunctions of an open packing are admissible, so $G_f\le G(O)\le0$. Delete every component with $\int|\nabla u_i|^2\ge\Lambda$. For a remaining component, functional Faber--Krahn gives
\[
 b\le |\{u_i>0\}|\int|\nabla u_i|^2
 <\Lambda|\{u_i>0\}|,
 \qquad
 n\le\frac{\Lambda|O|}{b},\qquad
 \mathcal E(U)\ge-\frac{\Lambda^2|O|}{b}.
\]
The positivity measures here are independent of the choice of Sobolev representatives. If $G_f=0$, then $G(O)=0$ and the empty packing proves the assertion. Suppose henceforth that $G_f<0$. A minimizing sequence can be trimmed as above and, after a subsequence, has one fixed count $n\ge1$ and $\sum_i\|u_i\|_{H^1}^2\le n(1+\Lambda)$. Boundedness of $O$ and Rellich compactness give, for $1\le i\le n$,
\[
 u_i^{(\nu)}\rightharpoonup u_i\quad\text{in }H_0^1(O),\qquad
 u_i^{(\nu)}\longrightarrow u_i\quad\text{in }L^2(O)\text{ and a.e.}
\]
Thus $\|u_i\|_2=1$, $u_i\ge0$, and $u_i u_k=0$ for $i\ne k$, while
\[
 G_f\le\mathcal E(U)\le\liminf_{\nu\to\infty}\mathcal E(U^{(\nu)})=G_f.
\]
Choose an attaining family with the fewest components and set $\mu_i=\int|\nabla u_i|^2$. Deletion gives
\[
 0\le\mathcal E(U\setminus\{u_i\})-\mathcal E(U)=\Lambda-\mu_i.
\]
Equality would give an attaining family with fewer components; hence $\mu_i<\Lambda$ for every $i$.

Write $\mathcal R(v)=\int|\nabla v|^2/\int v^2$ for $v\ne0$. For real $\phi\in C_c^\infty(O)$, replacing $u_i$ by $(1+t\phi)u_i/\|(1+t\phi)u_i\|_2$ is admissible for both signs of sufficiently small $t$. Consequently,
\[
 0=\frac12\left.\frac{d}{dt}\right|_{t=0}\mathcal R((1+t\phi)u_i)
 =\int\nabla u_i\cdot\nabla(\phi u_i)-\mu_i\int\phi u_i^2,
\]
which is \eqref{fullsquare:mult}. With $\rho,g,P$ as in Lemma~\ref{fullsquare:gradient}, summation yields
\[
 \Delta\rho=2g-2P\ge-2\Lambda\rho
 \quad\text{in }\mathcal D'(O).
\]
The function $e^{\sqrt{2\Lambda}z}\rho(x)$ is nonnegative and distributionally subharmonic in three dimensions. Its mean-value bound gives
\[
 \|\rho\|_{L^\infty(B_r(x))}
 \le C_{\Lambda,r}\|\rho\|_{L^1(B_{2r}(x))}<\infty
 \qquad(B_{2r}(x)\Subset O).
\]

For $X\in C_c^\infty(O;\mathbb R^2)$, let $\Phi_t$ be its flow and put
\[
 v_{i,t}=u_i\circ\Phi_t^{-1},\qquad
 A_i(t)=\int_O|\nabla v_{i,t}|^2,\qquad
 M_i(t)=\int_Ov_{i,t}^2.
\]
The family $(v_{i,t}/\sqrt{M_i(t)})_i$ is admissible. Change of variables gives
\begin{align*}
 A_i'(0)&=\int_O\left[|\nabla u_i|^2\operatorname{div}X
                   -2\nabla u_i\cdot DX\,\nabla u_i\right],\\
 M_i'(0)&=\int_Ou_i^2\operatorname{div}X.
\end{align*}
Since $A_i(0)=\mu_i$ and $M_i(0)=1$, minimality implies
\begin{equation}\label{fullsquare:inner-variation}
 \begin{aligned}
 0&=\left.\frac{d}{dt}\right|_{t=0}\sum_i\frac{A_i(t)}{M_i(t)}\\
  &=\int_O\left[(g-P)\operatorname{div}X
          -2\sum_i\nabla u_i\cdot DX\,\nabla u_i\right]
   =\int_O T:DX.
 \end{aligned}
\end{equation}
Thus $\operatorname{div}T=0$. Lemma~\ref{fullsquare:gradient} now gives locally Lipschitz representatives of all components.

Set $D_i=D_O(u_i)$ and use these representatives throughout. The sets $D_i$ are nonempty, open, and disjoint, and Lemma~\ref{fullsquare:positivity-trace} gives $u_i\in H_0^1(D_i)$ and the corresponding restriction property on every connected component. For $\psi\in C_c^\infty(D_i)$, both signs of the normalized variation $u_i+t\psi$ are admissible for small $|t|$, since $u_i$ has a positive minimum on $\operatorname{supp}\psi$. Hence
\[
 \int\nabla u_i\cdot\nabla\psi=\mu_i\int u_i\psi,
 \qquad -\Delta u_i=\mu_i u_i\quad\text{in }D_i.
\]
On every connected component $C$ of $D_i$, the nonzero function $u_i\mathbf1_C\in H_0^1(C)$ is a positive eigenfunction, so $\lambda_1(C)=\mu_i$. If $C_1,C_2$ were two such components, replacing $u_i$ by their normalized restrictions would change the energy by
\[
 2(\mu_i-\Lambda)-(\mu_i-\Lambda)=\mu_i-\Lambda<0.
\]
This contradicts minimality. Thus $D_i$ is connected, $\lambda_1(D_i)=\mu_i$, and
\[
 G(O)\le\sum_i(\lambda_1(D_i)-\Lambda)=G_f\le G(O).
\]

Finally, suppose that $O$ is connected and $G(O)<0$. The family is nonzero. For every center $x\in O$, radial testing of \eqref{fullsquare:inner-variation}, using $\operatorname{tr}T=-2P$, gives for almost every admissible radius
\[
 r\int_{\partial B_r(x)}
       \left(g-P-2\sum_i|\partial_r u_i|^2\right)
 =-2\int_{B_r(x)}P.
\]
The function $r\mapsto\int_{B_r(x)}g$ is absolutely continuous, and therefore
\[
 \frac{d}{dr}\int_{B_r(x)}g
 =2\int_{\partial B_r(x)}\sum_i|\partial_r u_i|^2
  +\int_{\partial B_r(x)}P-\frac2r\int_{B_r(x)}P.
\]
This is the Pohozaev law in \cite[Definition~1.2]{TavaresTerracini2012} for $f_i(s)=\mu_i s$: indeed,
\[
 \frac2r\int_{B_r(x)}\sum_i\mu_i u_i\nabla u_i\cdot(y-x)
 =\int_{\partial B_r(x)}P-\frac2r\int_{B_r(x)}P.
\]
Choose a positive-phase point $x_0$ and a connected exhaustion $x_0\in O_\ell\Subset O$, $\bigcup_\ell O_\ell=O$. On each $O_\ell$, the family is finite, nonzero, nonnegative, segregated, and Lipschitz, satisfies the positive-phase equations, and has the displayed Pohozaev law. Thus \cite[Theorem~1.1]{TavaresTerracini2012}, with $Z=\{U=0\}$, gives
\[
 \dim_{\mathrm H}(Z\cap O_\ell)\le1,
 \qquad
 \left|O\setminus\bigcup_iD_i\right|
 =|Z|\le\sum_{\ell=1}^\infty|Z\cap O_\ell|=0.
\]
The regularity theorem is applied only to a fixed finite vector; none of its quantitative constants is used for sequences with varying phase counts.
\end{proof}

\begin{lemma}[Coupling inequalities]\label{fullsquare:S}
An exact minimizing function family satisfies
\begin{equation}\label{fullsquare:Sformula}
 \begin{gathered}
 -\Delta u_i\le\mu_i u_i,\qquad
 -\Delta(2u_i-v)\ge2\mu_i u_i-w,\\
 v=\sum_k u_k=|U|,\qquad w=\sum_k\mu_k u_k.
 \end{gathered}
\end{equation}
\end{lemma}
\begin{proof}
The first inequality holds for any zero-extended nonnegative first eigenfunction. Test its equation with $\phi\min(u_i/s,1)$ for $\phi\ge0$, discard the nonnegative term, and let $s\downarrow0$.

For the second inequality, take a nonnegative smooth compactly supported $\phi$ and set
\[
 v_i=(2u_i-v+t\phi)_+,\qquad v_k=(u_k-t\phi)_+\quad(k\ne i).
\]
This is an admissible segregated $H_0^1$ family, and for a fixed finite family all its masses are positive when $t$ is small. Pointwise expansion gives
\[
 \sum_k\int|\nabla v_k|^2-\sum_k\mu_k
 =2t\int\nabla(2u_i-v)\cdot\nabla\phi+t^2\int|\nabla\phi|^2.
\]
The mass increments are $2t\int u_i\phi+r_i$ and $-2t\int u_k\phi+r_k$, respectively, where $r_k\ge0$ and $\sum r_k=t^2\int\phi^2$. Continuity of the positive-part map in $H^1$ makes each energy increment $o(1)$. Expanding the normalization denominators to first order therefore has total error $o(t)$. Exact minimality yields
\[
 0\le2t\left\{\int\nabla(2u_i-v)\cdot\nabla\phi
 -\int(2\mu_i u_i-w)\phi\right\}+o(t).
\]
Let $t\downarrow0$ with the finite family fixed. The resulting distributional inequality is exact and can subsequently be used for varying counts without any uniformity requirement on the discarded $o(t)$.
\end{proof}

\subsection{Uniform estimates for reflected families}
Take an exact minimizing family in $Q_R$ from Lemma~\ref{fullsquare:existence}, using its canonical representatives and positivity sets. For large $R$, a disk with first eigenvalue below $\Lambda$ shows that $G(Q_R)<0$. Reflect across the square edges to tile the plane, assigning independent labels to all mirror copies and preserving their unit masses. For fixed $R$, every compact set meets only finitely many labels.

\begin{lemma}[Bounds for the reflected configuration]\label{fullsquare:reflection}
The reflected family satisfies \eqref{fullsquare:Sformula}, \eqref{fullsquare:mult}, and $\operatorname{div}T=0$ on the whole plane. There exist a threshold $R_{\rm cmp}=R_{\rm cmp}(\Lambda)>0$ and constants $K,L,C_0,\eta_0>0$, depending only on $\Lambda$, such that, for all $R\ge R_{\rm cmp}(\Lambda)$ and all centers $x$,
\begin{equation}\label{fullsquare:uniform}
 \rho\le K^2,\qquad g\le L^2,\qquad
 \int_{B_{C_0}(x)}\rho\ge\eta_0,\qquad
 \Delta g+2\Lambda g\ge0.
\end{equation}
\end{lemma}
\begin{proof}
Each reflected component is a complete zero-extended first eigenfunction, so \eqref{fullsquare:mult} holds on the plane. Set
\[
 \sigma_i=\mu_i u_i+\Delta u_i\ge0.
\]
The zero-extension inequality proved in Lemma~\ref{fullsquare:S} makes these nonnegative Radon measures in $H^{-1}_{\rm loc}$. The coupling inequality inside the base square transfers by reflection, so on the interior of each reflected square
\[
 \sum_k\sigma_k-2\sigma_i\ge0.
\]
Let $e$ be an open grid edge and let $\mathcal R_e$ denote reflection across its line. For a label $i$ in either adjacent square, let $i^*\ne i$ be its independent mirror label. Reflection of the distributional equation gives
\[
 \sigma_{i^*}=(\mathcal R_e)_\#\sigma_i,
 \qquad \sigma_{i^*}|_e=\sigma_i|_e,
 \qquad
 \left(\sum_k\sigma_k-2\sigma_i\right)\Big|_e\ge0.
\]
For labels in other squares the restriction to $e$ is zero. The sums are locally finite. A point has zero $H^1$ capacity in two dimensions; consequently a nonnegative $H^{-1}_{\rm loc}$ measure has no atom at a grid vertex. The square interiors, open edges, and vertices exhaust the plane, so $\sum_k\sigma_k\ge2\sigma_i$ globally. Together with $\sigma_i\ge0$, this is \eqref{fullsquare:Sformula}.

Let $T_0$ denote the stress tensor of the minimizing configuration on the base square $Q_R$, defined by the formula in Lemma~\ref{fullsquare:gradient}. To verify stress conservation, write $Q_a=F_a(Q_R)$ for the reflected squares and $S_a=DF_a$. For $X\in C_c^\infty(\mathbb R^2;\mathbb R^2)$, define on the base square
\[
 Y=\sum_a S_a^{\mathsf T}X\circ F_a,
 \qquad
 DY=\sum_a S_a^{\mathsf T}(DX\circ F_a)S_a.
\]
Only finitely many terms occur. On a base edge with normal $\nu$, the terms pair so that $F_a x=F_{a'}x$ and $S_{a'}\nu=-S_a\nu$; hence $Y\cdot\nu=0$. At a corner both normal components vanish. Thus the flow of $Y$ preserves $Q_R$. The change-of-variables calculation in \eqref{fullsquare:inner-variation} applies to this flow and gives $\int_{Q_R}T_0:DY=0$. Since $T(F_a x)=S_aT_0(x)S_a^{\mathsf T}$ almost everywhere,
\[
 \int_{\mathbb R^2}T:DX
 =\sum_a\int_{Q_R}T_0:
       S_a^{\mathsf T}(DX\circ F_a)S_a
 =\int_{Q_R}T_0:DY=0.
\]

Lemma~\ref{holder:heat-density} follows from Dirichlet heat-semigroup domination on arbitrary open sets in Appendix~\ref{app:heat}; its proof is independent of the compactness and positive lower bounds in Sections~\ref{sec:compactness}--\ref{fullsquare:sec}. It gives
\[
 0\le u_i\le K:=\sqrt{e\Lambda/(4\pi)},
 \qquad \rho=\sum_i u_i^2\le K^2.
\]
The last inequality uses segregation and has no factor counting labels. Testing \eqref{fullsquare:mult} with $\eta^2$ and summing yields
\[
 \int\eta^2g\le2\Lambda\int\eta^2\rho
                    +4\int|\nabla\eta|^2\rho.
\]
For cutoffs on fixed-radius balls this is uniform in the center and in $R$. Lemma~\ref{fullsquare:gradient}, applied on every ball $B_2(x)$, now gives
\[
 g\le L^2:=C(\Lambda)K^2
 \quad\text{a.e. on }\mathbb R^2.
\]
In particular, the reflected functions have locally Lipschitz representatives across the grid edges as well.

To prove the distributional inequality for $g$, fix one reflected family. On each member of a bounded connected exhaustion containing a fixed nonzero point, only finitely many labels occur. The positive-phase equations and radial stress identity verify the assumptions of \cite[Theorem~1.1]{TavaresTerracini2012}, as in Lemma~\ref{fullsquare:existence}. The common zero set therefore consists of regular arcs and a relatively closed, locally finite singular set $\mathcal S$.

The two sides of a regular arc have distinct labels. Otherwise, in a small neighborhood only one label $u_i$ is present, and the two coupling inequalities reduce to
\[
 -\Delta u_i\le\mu_i u_i,
 \qquad -\Delta u_i\ge\mu_i u_i.
\]
The strong maximum principle contradicts the interior zero of this nonnegative, nonzero solution. For two distinct labels $i,k$, applying the second coupling inequality to both labels gives, with $v=u_i-u_k$,
\[
 -\Delta v=\mu_i v^+-\mu_kv^-,\qquad g=|\nabla v|^2.
\]
Interior elliptic regularity and the Sobolev chain rule give the weak Bochner identity
\[
 \Delta g
 =2|D^2v|^2
  -2\bigl(\mu_i\mathbf1_{\{v>0\}}+\mu_k\mathbf1_{\{v<0\}}\bigr)|\nabla v|^2
 \ge-2\Lambda g.
\]
The same inequality holds inside each positive phase, and hence on $\mathbb R^2\setminus\mathcal S$.

On a relatively compact disk $D$, choose a Newton potential $w$ with $\Delta w=2\Lambda g$. Since $g$ is bounded, $w$ is continuous and locally bounded. The function $h=g+w$ is subharmonic on $D\setminus\mathcal S$ and locally bounded above. The locally finite set $\mathcal S$ is polar; see \cite[Example~5.1.2, p.~124]{ArmitageGardiner2001}. Applying \cite[Corollary~5.2.2, p.~128]{ArmitageGardiner2001} to $-h$ gives the subharmonic extension
\[
 \widetilde h(x)=\limsup_{\substack{y\to x\\y\notin\mathcal S}}h(y),
 \qquad \widetilde h=h\quad\text{a.e. on }D.
\]
Therefore, for every nonnegative $\phi\in C_c^\infty(D)$,
\[
 \langle\Delta g+2\Lambda g,\phi\rangle
 =\int_D(g+w)\Delta\phi
 =\int_D\widetilde h\,\Delta\phi\ge0.
\]
This proves the inequality on the whole plane.

For the mass bound, fix $r=r(\Lambda)>0$ such that $\lambda_1(B_r)<\Lambda/2$. For $B_{2r}(y)\Subset Q_R$, choose a smooth cutoff with
\[
 0\le\chi\le1,\qquad
 \chi=0\text{ on }B_r(y),\qquad
 \chi=1\text{ outside }B_{2r}(y),\qquad
 |\nabla\chi|^2\le C/r^2.
\]
Set $a=\int_{B_{2r}(y)}\rho$. Suppose first that $a<1/2$. For every base label, let $a_i=\int(1-\chi^2)u_i^2\le a$. Testing its eigenfunction equation with $\chi^2u_i$ gives
\[
 \int|\nabla(\chi u_i)|^2
 =\mu_i\int\chi^2u_i^2+\int|\nabla\chi|^2u_i^2,
 \qquad \int\chi^2u_i^2=1-a_i\ge1-a.
\]
Writing $\mathcal R(v)=\int|\nabla v|^2/\int v^2$, we obtain
\begin{equation}\label{fullsquare:mass-cutoff-cost}
 \sum_i\bigl(\mathcal R(\chi u_i)-\mu_i\bigr)
 =\sum_i\frac{\int|\nabla\chi|^2u_i^2}{1-a_i}
 \le\frac{Ca}{r^2(1-a)}.
\end{equation}
For each $i$, select a connected component $W_i$ of $\{\chi u_i>0\}$ for which
\[
 \lambda_1(W_i)\le\mathcal R(\chi u_i\mathbf1_{W_i})
 \le\mathcal R(\chi u_i).
\]
Such a component exists by the mass-weighted decomposition of the quotient; Lemma~\ref{fullsquare:positivity-trace} supplies its zero trace. The sets $W_i$ form a packing disjoint from $B_r(y)$. Inserting the disk and using exact minimality gives
\begin{align*}
 0&\le\sum_i\bigl(\lambda_1(W_i)-\mu_i\bigr)
               +\lambda_1(B_r)-\Lambda\\
  &\le\frac{Ca}{r^2(1-a)}-\frac\Lambda2.
\end{align*}
Consequently, in both cases $a<1/2$ and $a\ge1/2$,
\[
 \int_{B_{2r}(y)}\rho\ge
 \eta:=\min\left\{\frac12,\frac{\Lambda r^2}{2C+\Lambda r^2}\right\}>0.
\]
The constants depend only on $\Lambda$ and the fixed cutoff profile.

For $R\ge6r$ and any $x\in\mathbb R^2$, choose a grid square whose closure contains $x$. In its square coordinates, project each coordinate of $x$ onto $[3r,R-3r]$ to obtain $y$. Then
\[
 |x-y|\le3\sqrt2\,r,\qquad
 B_{2r}(y)\Subset Q_a,\qquad
 B_{2r}(y)\subset B_{C_0}(x),\qquad C_0=(2+3\sqrt2)r.
\]
The preceding estimate transfers to $Q_a$ by reflection. Thus $\int_{B_{C_0}(x)}\rho\ge\eta$ for every $x$, and one may take $\eta_0=\eta$ and $R_{\rm cmp}(\Lambda)=6r$.
\end{proof}

\subsection{Strong compactness under a single relabeling}
As the square grows, coordinatewise compactness alone does not control the total energy. We first bound the contribution of labels of small mass, then recover the limit energy identity from the coupling inequalities.

\begin{lemma}[Compactness with varying label count]\label{fullsquare:compactness}
Consider reflected families from squares with side lengths $R_j\ge R_{\rm cmp}(\Lambda)$, where $R_{\rm cmp}(\Lambda)$ is the threshold in Lemma~\ref{fullsquare:reflection}. After arbitrary translations, any sequence of these families has a subsequence that converges strongly in $H^1_{\rm loc}(\ell^2)$ and $C^0_{\rm loc}(\ell^2)$ after one permutation of the labels in each member. The limit is nonzero and satisfies the coupling inequalities, multiplication identity, stress conservation, and the corresponding local bounds in \eqref{fullsquare:uniform}. Each limit label belongs globally to $H^1(\mathbb R^2)$ and has mass at most one.
\end{lemma}
\begin{proof}
Order the labels by nonincreasing values of $s_{j,i}=\int e^{-|x|}u_{j,i}^2$. The density bound gives $\sum_i s_{j,i}\le S$, hence $s_{j,i}\le S/(N+1)$ for $i>N$. Take a smooth cutoff $0\le\psi\le1$, equal to one on $B_{2r}$ and supported in $B_{4r}$, and set
\[
 p_{j,i}=\int\psi^2u_{j,i}^2,\qquad
 e_{j,i}=\int|\nabla(\psi u_{j,i})|^2.
\]
Functional Faber--Krahn and segregation imply
\[
 \sum_i\sqrt{p_{j,i}}
 \le\left(\frac{|B_{4r}|}{b}\sum_i e_{j,i}\right)^{1/2}\le C_r.
\]
Each tail label has $p_{j,i}\le e^{4r}S/(N+1)$. Consequently,
\begin{equation}\label{fullsquare:mass-tail}
 \sup_j\sum_{i>N}p_{j,i}
 \le \sup_j\left[
   \left(\sup_{i>N}\sqrt{p_{j,i}}\right)
   \sum_i\sqrt{p_{j,i}}\right]
 \le\frac{C_r}{\sqrt{N+1}}.
\end{equation}
Choose a smooth cutoff $\eta=1$ on $B_r$, supported in $B_{2r}$, with $|\nabla\eta|\le C/r$. The multiplication identity and Young's inequality give
\[
 \int\eta^2|\nabla u_{j,i}|^2
 \le 2\Lambda\int\eta^2u_{j,i}^2
       +4\int|\nabla\eta|^2u_{j,i}^2
 \le C_{r,\Lambda}p_{j,i}.
\]
Together with \eqref{fullsquare:mass-tail}, this proves the uniform tail estimate
\begin{equation}\label{fullsquare:energy-tail}
 \sup_j\sum_{i>N}\|u_{j,i}\|_{H^1(B_r)}^2
 \le\frac{C_r}{\sqrt{N+1}}.
\end{equation}
Here and below $C_r$ is independent of $j,N$. Coordinatewise Rellich compactness, a diagonal selection, and \eqref{fullsquare:mass-tail} yield strong local $L^2(\ell^2)$ and weak local $H^1(\ell^2)$ convergence. After a further subsequence, $\mu_{j,i}\to\mu_i$ for each $i$. Weak lower semicontinuity on finite blocks, followed by monotone convergence in the block size, gives the same tail bound for the limit. Each original label has global $H^1$ norm at most $\sqrt{1+\Lambda}$ and unit mass, so each limit label has global mass at most one.

The gradient of the total amplitude $v_j=|U_j|$ is the sum of the segregated label gradients, and its tail satisfies the same bound. Thus the coupling inequalities pass to the limit. To upgrade gradient convergence, we derive the energy identity in the limit. Both sides of the distributional coupling inequalities act continuously on compactly supported $H_0^1$ test functions; nonnegative smooth approximation extends the inequalities to nonnegative tests of this type. Test with $\eta^2u_i$, where $\eta$ is nonnegative, smooth, and compactly supported. Segregation and vanishing of Sobolev gradients on zero sets give
\[
 \nabla v\cdot\nabla(\eta^2u_i)
 =\nabla u_i\cdot\nabla(\eta^2u_i),\qquad
 wu_i=\mu_i u_i^2\quad\text{a.e.}
\]
The two coupling inequalities give opposite inequalities for the same expression. Hence
\begin{equation}\label{fullsquare:limit-energy}
 \int\eta^2|\nabla u_i|^2
 =\mu_i\int\eta^2u_i^2
 -2\int\eta u_i\nabla\eta\cdot\nabla u_i.
\end{equation}
The same test with an arbitrary nonnegative smooth weight in place of $\eta^2$, followed by taking differences, gives the full multiplication identity \eqref{fullsquare:mult}.

The original sequence satisfies \eqref{fullsquare:limit-energy} as well. Convergence of $\mu_{j,i}$, strong local $L^2$ convergence, and weak gradient convergence pass every term on the right to the limit. Thus the weighted gradient norms converge, and every fixed finite block converges strongly in local $H^1$. For each fixed $r,N$, \eqref{fullsquare:energy-tail} and its limit counterpart now give
\begin{align*}
 \limsup_{j\to\infty}\|U_j-U\|_{H^1(B_r;\ell^2)}^2
 &\le \limsup_{j\to\infty}
       \sum_{i\le N}\|u_{j,i}-u_i\|_{H^1(B_r)}^2\\
 &\quad+2\limsup_{j\to\infty}
       \sum_{i>N}\bigl(\|u_{j,i}\|_{H^1(B_r)}^2
                       +\|u_i\|_{H^1(B_r)}^2\bigr)\\
 &\le \frac{C_r}{\sqrt{N+1}}.
\end{align*}
First $j\to\infty$ and then $N\to\infty$ therefore give strong local $H^1(\ell^2)$ convergence. The uniform Lipschitz bound gives local uniform convergence: a vector difference of size $a$ at one point forces a fixed positive $L^2$ difference on a ball of radius $a/(4L)$. Stress conservation and the distributional inequality for $g$ pass by strong gradient convergence. The fixed-ball lower mass bound at every center ensures that the limit is nonzero.
\end{proof}

\subsection{Frequency and normalized blow-ups for a countable limit}
At this stage the limit may still have countably many labels. Frequency controls its vanishing order before any local finiteness has been established. The half-integer value will be used in the complex-stress expansion below.

\begin{lemma}[Frequency of the limit]\label{fullsquare:frequency}
For any center of a nonzero limit as above, write
\[
 H(r)=r^{-1}\int_{\partial B_r}\rho,\qquad
 D(r)=\int_{B_r}(g-P),\qquad F(r)=\int_{B_r}P.
\]
Then $H(r)>0$ for $0<r<\Lambda^{-1/2}$, and $N(r)=D(r)/H(r)$ satisfies
\begin{equation}\label{fullsquare:frequency-formula}
 H'=2NH/r,\qquad N\ge-\Lambda r^2,\qquad
 \frac{N+1}{1-\Lambda r^2}\text{ is nondecreasing}.
\end{equation}
At a common zero, $\alpha=N(0+)=q/2$ for an integer $q\ge2$.
\end{lemma}
\begin{proof}
The multiplication identity and radial stress tests give, for almost every radius,
\[
 D=\int_{\partial B_r}\sum_i u_i\partial_ru_i,
 \qquad H'=2D/r,\qquad
 rD'=2r\int_{\partial B_r}|\partial_rU|^2-2F.
\]
Consequently,
\[
 N'=\frac2H\left(\int_{\partial B_r}|\partial_rU|^2-
 \frac{D^2}{rH}\right)-\frac{2F}{rH}.
\]
Let $I=\int_{B_r}\rho$. The divergence identity and the Cauchy--Young inequality give
\[
 2I=r^2H-2\int_{B_r}U\cdot(x\cdot\nabla U)
 \le r^2H+I+r^2(D+F).
\]
Hence $I\le r^2(H+D+F)$, while $F\le\Lambda I$. Thus
\[
 \frac FH\le\frac{\Lambda r^2(1+N)}{1-\Lambda r^2}.
\]
Together with $D\ge-F$, this implies $N\ge-\Lambda r^2$ and \eqref{fullsquare:frequency-formula}.

Consider a connected component $(a,b)$ of $\{r\in(0,\Lambda^{-1/2}):H(r)>0\}$, and fix $s$ in it. Monotonicity gives $N(r)\le C_s$ for $r\le s$, hence $H(r)\ge H(s)(r/s)^{2C_s}$. If $a>0$, continuity extends positivity through $a$, a contradiction. Thus $a=0$. For $r\ge s$, the lower frequency bound gives $H(r)\ge H(s)e^{-\Lambda(r^2-s^2)}$, so the same argument gives $b=\Lambda^{-1/2}$.

If the interior $Z$ of the common zero set were nonempty, then at every center in $Z$, $H$ would vanish for small radii. The component argument forces it to vanish throughout $(0,\Lambda^{-1/2})$. The entire ball of that fixed radius would then lie in $Z$. This fixed-radius property makes $Z$ both open and closed in the plane, contrary to the nonzero limit. Hence $H>0$ as claimed. Monotonicity now yields a finite limit $\alpha=N(0+)\ge0$. At a common zero, Lipschitz continuity gives $H(r)\le Cr^2$; integrating the frequency identity implies $\alpha\ge1$.

We prove half-integrality without assuming finitely many local labels. For any $r_j\downarrow0$, define
\[
 V_j(x)=\frac{U(x_0+r_jx)}{\sqrt{H(r_j)}}.
\]
Bounded frequency and $H'/H=2N/r$ give uniform total $H^1$ bounds on every fixed ball. The source coefficients $r_j^2\mu_i$ tend to zero uniformly. Choose $0\le\eta\le1$ supported in $B_8$, equal to one on $B_4$, and put $p_{j,i}=\|\eta v_{j,i}\|_2^2$. Functional Faber--Krahn and the energy bound give $\sum_i\sqrt{p_{j,i}}\le C$. Discard zero components when restricting to $B_4$, and order the positive $p_{j,i}$ once in each member. Then
\[
 \sqrt{p_{j,N+1}}\le\frac C{N+1},\qquad
 \sum_{i>N}p_{j,i}\le\frac{C^2}{N+1}.
\]
Multiplication cutoffs give the same gradient-tail order on $B_3$. The coupling and energy-identity argument of Lemma~\ref{fullsquare:compactness} yields a strong $H^1(B_3;\ell^2)$ sublimit $V_0$. This relabeling is used only on these fixed nested neighborhoods. The rescaled inequality $\Delta g_j+2\Lambda r_j^2g_j\ge0$, combined with the energy bound and the three-dimensional subharmonic mean-value estimate, gives an interior uniform gradient bound. Hence convergence is also locally uniform on $B_2$. The trace of strong $H^1$ convergence gives $H_0(1)=1$, so the limit is nonzero.

For almost every fixed $0<t<2$, strong convergence and the frequency ratios give $H_0(t)=t^{2\alpha}$ and $N_0(t)=\alpha$. The zero-potential frequency identity is therefore
\[
 0=N_0'(t)=\frac{2}{H_0(t)}
 \int_{\partial B_t}
 \left|\partial_rV_0-\frac{\alpha}{t}V_0\right|_{\ell^2}^2.
\]
Thus $\partial_r(r^{-\alpha}V_0)=0$ almost everywhere, and continuity gives
\[
 V_0(r,\theta)=r^\alpha\Phi(\theta).
\]
On each positive angular interval, the corresponding component is harmonic, so $\phi_i''+\alpha^2\phi_i=0$. A positive interval has zero values at both endpoints and length exactly $\pi/\alpha$. A component cannot be positive on the whole circle, since integrating the differential equation would contradict $\alpha>0$. There are therefore only finitely many positive angular intervals. The same exclusion of open common zero sets for the zero-potential frequency shows that no open angular interval is left uncovered. Thus $q\pi/\alpha=2\pi$, or $\alpha=q/2$. The homogeneous expression, first obtained on a fixed ball, extends homogeneously while preserving the local equations and stationarity, so this last argument does not require a common global relabeling of the blow-up sequence.
\end{proof}

\subsection{Complex stress, strict label exclusion, and cell diameters}
Local seam comparisons must remove complete cells. We first show that all labels near a compact set eventually come from a fixed finite list. A uniform lower mass bound for each visible label then bounds complete-cell diameters.

\begin{lemma}[Strict local label exclusion]\label{fullsquare:clearing}
Let $U_j\to U$ be a translated sequence and its strong limit as in Lemma~\ref{fullsquare:compactness}. For every compact set $K$ and open neighborhood $V$ of $K$, there are a finite set $J$ of nonzero limit labels and an index $j_0$ such that
\[
 u_{j,i}|_K=u_i|_K=0\quad(i\notin J,\ j\ge j_0),
 \qquad \sup_V u_i>0\quad(i\in J).
\]
In particular, near each compact set only finitely many limit labels occur, and the approximating labels there eventually belong to a fixed finite set of limit labels nonzero in $V$.
\end{lemma}
\begin{proof}
For the sequence and its limit, use the complex stress
\[
 \mathcal Q=\sum_i(u_{i,x}-\mathrm i u_{i,y})^2,
 \qquad |\mathcal Q|=g,\qquad
 \bar\partial\mathcal Q=-\partial P,
 \qquad |\nabla P|\le2\Lambda|U|\sqrt g.
\]
Fix $p>2$. The uniform bounds for $U_j$ and $\nabla U_j$, together with the local Cauchy--Riemann estimates used in Lemma~\ref{fullsquare:gradient}, give uniform $W^{1,p}$ bounds for $\mathcal Q_j$ on every compactly contained disk. Moreover,
\[
 \|\mathcal Q_j-\mathcal Q\|_{L^1(B)}
 \le\bigl(\|\nabla U_j\|_{L^2(B)}+\|\nabla U\|_{L^2(B)}\bigr)
       \|\nabla U_j-\nabla U\|_{L^2(B)}\longrightarrow0.
\]
Morrey compactness and uniqueness of the $L^1$ limit therefore imply $\mathcal Q_j\to\mathcal Q$ locally uniformly.

For each $x_*\in K$, we will find an open neighborhood $W_*\Subset V$, a finite label set $J_*$, and $j_*$ such that
\begin{equation}\label{fullsquare:local-exclusion}
 x_*\in W_*,\qquad
 u_{j,i}|_{W_*}=0\quad(i\notin J_*,\ j\ge j_*),
\end{equation}
with a fixed point $z_i\in V$ satisfying $u_i(z_i)>0$ for each $i\in J_*$. If $u_\ell(x_*)>0$, choose $W_*\Subset V$ and $\gamma>0$ with $u_\ell\ge2\gamma$ on $\overline W_*$. Local uniform convergence gives $u_{j,\ell}\ge\gamma$ there for all large $j$; segregation proves \eqref{fullsquare:local-exclusion} with $J_*=\{\ell\}$ and $z_\ell=x_*$.

\emph{Nondegenerate common zeros.}
Suppose that $U(x_*)=0$ and $\mathcal Q(x_*)\ne0$. Translate $x_*$ to the origin. On a sufficiently small disk compactly contained in $V$, choose consistent square roots
\[
 p_j=\sqrt{\mathcal Q_j}\longrightarrow p=\sqrt{\mathcal Q},
 \qquad b_j=(\operatorname{Re}p_j,-\operatorname{Im}p_j),
 \qquad b=(\operatorname{Re}p,-\operatorname{Im}p).
\]
The convergence is uniform, and the roots are uniformly bounded in $W^{1,p}$, since the stresses stay away from zero. After a fixed rotation and a choice of sign, shrink to a square $D=(-d,d)^2\Subset V$ on which, for all sufficiently large $j$ and for the limit,
\[
 (b_j)_1\ge\beta>0,\qquad |b_j|\le M,
 \qquad |\mathcal Q_j|\ge c>0.
\]
On a connected positive piece, the phase equation and the identity defining $\mathcal Q_j$ give
\[
 \nabla u_{j,i}=\epsilon b_j,\qquad \epsilon\in\{-1,1\},
\]
with a constant sign on that piece: the eigenfunction gradient is continuous and $b_j$ is continuous and nonzero. For the limit, the two coupling inequalities give the phase equation on $\{u_i>0\}$, so the same argument applies. Gradients vanish almost everywhere on Sobolev zero sets. Hence the common zero set in $D$ has area zero, since $g_j=|\mathcal Q_j|\ge c$, and
\[
 \operatorname{curl}b_j=0\quad\text{a.e. in }D.
\]
Because $b_j\in W^{1,p}(D)$, its distributional curl has no additional singular part. The Poincare lemma on the square provides potentials normalized by
\[
 \nabla s_j=b_j,\qquad s_j(0)=0,
 \qquad s_j\longrightarrow s\quad\text{in }C^1(\overline D).
\]
Here $\overline D$ is chosen inside the disk on which the roots were defined.

Set
\[
 a=\frac{\beta d}{2},\qquad
 b_0=\min\left\{\frac d2,\frac{\beta d}{4M}\right\},\qquad
 \mathcal R=(-a,a)\times(-b_0,b_0),\qquad
 F_j(x,y)=(s_j(x,y),y).
\]
Since $\partial_xs_j\ge\beta$, each $F_j$ is one-to-one on $D$. For $|y|\le b_0$,
\[
 s_j(d,y)\ge\beta d-Mb_0\ge\frac{3\beta d}{4}>a,
 \qquad
 s_j(-d,y)\le-\frac{3\beta d}{4}<-a.
\]
Thus $\overline{\mathcal R}\subset F_j(D)$ for all large $j$, and also $\overline{\mathcal R}\subset F(D)$, where $F(x,y)=(s(x,y),y)$. These maps have $C^1$ inverses on the common rectangle.

For $\widehat u_{j,i}=u_{j,i}\circ F_j^{-1}$, the Sobolev chain rule gives, almost everywhere on $\mathcal R$,
\begin{align*}
 \partial_\sigma\widehat u_{j,i}
 &=\left(\frac{\partial_xu_{j,i}}{\partial_xs_j}\right)\circ F_j^{-1},\\
 \partial_\eta\widehat u_{j,i}
 &=\left(\partial_yu_{j,i}
       -\frac{\partial_ys_j}{\partial_xs_j}\partial_xu_{j,i}\right)
       \circ F_j^{-1}=0.
\end{align*}
Consequently $\widehat u_{j,i}(\sigma,\eta)=h_{j,i}(\sigma)$, and on each positive interval $I$,
\[
 h_{j,i}(\sigma)=\epsilon_I\sigma+c_I,
 \qquad \epsilon_I\in\{-1,1\}.
\]
An interval with both endpoints in $(-a,a)$ would have zero values at both endpoints, which is impossible for an affine function of slope $\pm1$. Every positive interval therefore reaches $-a$ or $a$. Segregation leaves at most two positive intervals across all labels. These conclusions hold for both the approximating families and the limit.

The common zero set has area zero, and $U(0)=0$. Hence the two positive intervals of the limit are exactly $(-a,0)$ and $(0,a)$, with respective labels $\ell_-,\ell_+$, possibly equal. In particular,
\[
 z_\pm=F^{-1}(\pm a/2,0)\in V,
 \qquad u_{\ell_\pm}(z_\pm)=a/2.
\]
Fix these points and the neighborhood
\[
 W_*=F^{-1}\bigl((-a/8,a/8)\times(-b_0/2,b_0/2)\bigr)\Subset D.
\]
Choose one index $j_*$ beyond the preceding thresholds such that, for every $j\ge j_*$,
\[
 \|s_j-s\|_{C^0(\overline D)}<a/8,
 \qquad
 \|U_j-U\|_{C^0(\overline D;\ell^2)}<a/4.
\]
Then
\[
 \begin{gathered}
 u_{j,\ell_\pm}(z_\pm)\ge a/4,\qquad
 s_j(z_-)<-3a/8,\qquad s_j(z_+)>3a/8,\\
 F_j(W_*)\subset(-a/4,a/4)\times(-b_0/2,b_0/2).
 \end{gathered}
\]
If $\ell_-\ne\ell_+$, the two possible labels are already present. If $\ell_-=\ell_+$, a different positive interval meeting $(-a/4,a/4)$ would be trapped strictly between the two witness coordinates and would have two interior zero endpoints. Both cases give
\[
 h_{j,i}|_{(-a/4,a/4)}=0
 \quad(i\notin\{\ell_-,\ell_+\},\ j\ge j_*).
\]
This proves \eqref{fullsquare:local-exclusion} with $J_*=\{\ell_-,\ell_+\}$ and fixed points $z_\pm\in V$.

\emph{Degenerate common zeros.}
Suppose that $U(x_*)=0$ and $\mathcal Q(x_*)=0$, and translate $x_*$ to the origin. By Lemma~\ref{fullsquare:frequency}, $\alpha=N(0+)=q/2$ with $q\ge2$, so $k=2\alpha-2$ is a nonnegative integer. Choose $r_0>0$ small enough that $\overline{B_{r_0}}\subset V$, $\Lambda r_0^2\le1$, and $|N(r)-\alpha|\le1/4$ for $0<r\le r_0$. Integration of $H'/H=2N/r$ gives
\begin{equation}\label{fullsquare:degenerate-growth}
 c_0r^{2\alpha+1/2}\le H(r)\le C_0r^{2\alpha-1/2}
 \qquad(0<r\le r_0).
\end{equation}
In particular,
\[
 \int_{B_r}\rho=\int_0^r tH(t)\,dt
 \le Cr^{2\alpha+3/2}.
\]
For $0<r<r_0/4$, the mean-value estimate for $\Delta g+2\Lambda g\ge0$ and the multiplication identity with a cutoff on $B_{4r}$ imply
\begin{align*}
 \sup_{B_r}|\mathcal Q|
 &=\operatorname*{ess\,sup}_{B_r}g
 \le Cr^{-2}\int_{B_{2r}}g\\
 &\le Cr^{-2}(r^{-2}+\Lambda)\int_{B_{4r}}\rho
 \le Cr^{2\alpha-5/2}=Cr^{k-1/2}.
\end{align*}
Since $U(0)=0$, the gradient bound on $B_r$ also gives
\[
 \sup_{B_r}|U|\le r\operatorname*{ess\,sup}_{B_r}\sqrt g
 \le Cr^{\alpha-1/4},
 \qquad
 |\bar\partial\mathcal Q(z)|\le C|z|^{k+1/2}
\]
for almost every sufficiently small $z$, using $\bar\partial\mathcal Q=-\partial P$.

Shrink $r_0$ once more so the last bound holds on $B_{r_0}$, and extend $f=\bar\partial\mathcal Q$ by zero outside this disk. Its Cauchy transform
\[
 Tf(z)=\pi^{-1}\int\frac{f(\zeta)}{z-\zeta}\,dA(\zeta)
\]
differs from $\mathcal Q$ by a holomorphic function in $B_{r_0}$. The coefficients
\[
 c_l=-\pi^{-1}\int f(\zeta)\zeta^{-l-1}\,dA(\zeta)
 \qquad(0\le l\le k+1)
\]
are absolutely integrable, since their radial majorants near zero are $C\rho^{k+1/2-l}$. For $r=|z|$, the inner part of the Taylor remainder is bounded by
\[
 \int_{|\zeta|<2r}\frac{|f(\zeta)|}{|z-\zeta|}\,dA
 +\sum_{l=0}^{k+1}r^l
   \int_{|\zeta|<2r}\frac{|f(\zeta)|}{|\zeta|^{l+1}}\,dA
 \le Cr^{k+3/2}.
\]
On $|\zeta|\ge2r$, the kernel identity
\[
 \frac1{z-\zeta}+\sum_{l=0}^{k+1}\frac{z^l}{\zeta^{l+1}}
 =\frac{z^{k+2}}{\zeta^{k+2}(z-\zeta)}
\]
bounds the outer remainder by
\[
 Cr^{k+2}\int_{2r}^{r_0}\rho^{-3/2}\,d\rho
 \le Cr^{k+3/2}.
\]
Adding the holomorphic Taylor expansion yields
\[
 \mathcal Q(z)=\sum_{l=0}^{k+1}a_lz^l+O(|z|^{k+3/2}).
\]
The bound $|\mathcal Q(z)|\le C|z|^{k-1/2}$ forces $a_l=0$ for $l<k$. If $a_k=0$, then $g\le C|z|^{k+1}$, and integration of the gradient from $U(0)=0$ gives
\[
 |U(z)|\le C|z|^{(k+3)/2},\qquad
 H(r)\le Cr^{k+3}=Cr^{2\alpha+1},
\]
contradicting \eqref{fullsquare:degenerate-growth} as $r\downarrow0$. Thus
\[
 \mathcal Q(z)=a_kz^k+O(|z|^{k+1}),\qquad a_k\ne0.
\]
The condition $\mathcal Q(0)=0$ implies $k\ge1$. For some fixed $r_1>0$,
\[
 |\mathcal Q(z)|\ge\tfrac12|a_k|\,|z|^k>0
 \qquad(0<|z|<r_1).
\]

Choose $r<\min\{r_1,j_{0,1}/\sqrt\Lambda\}$ with $\overline{B_{2r}}\subset V$. At each point of $\partial B_r$, either a component of $U$ is positive or $\mathcal Q$ is nonzero, so one of the preceding neighborhoods applies. Choose a finite subcover, let $J^\circ$ be the union of its finite label sets, and let $j^\circ$ be the maximum of its thresholds. For $j\ge j^\circ$,
\[
 u_{j,i}|_{\partial B_r}=0\qquad(i\notin J^\circ).
\]
If such a label had a cell $D_{j,i}$ meeting $B_{r/2}$, connectedness would force $D_{j,i}\subset B_r$: a path to any point outside $B_r$ would cross the circle at a positive point of that same label. Domain monotonicity would then give
\[
 \Lambda>\mu_{j,i}=\lambda_1(D_{j,i})
 \ge\lambda_1(B_r)=\frac{j_{0,1}^2}{r^2}>\Lambda,
\]
a contradiction. Hence \eqref{fullsquare:local-exclusion} holds with $W_*=B_{r/2}$, $J_*=J^\circ$, and $j_*=j^\circ$. All the chosen points lie in $V$ and were fixed before taking $j\to\infty$.

Finally, choose finitely many of the neighborhoods $W_*$ covering $K$, say $W_1,\ldots,W_M$, and set
\[
 J=\bigcup_{\ell=1}^M J_\ell,\qquad
 j_0=\max_{1\le\ell\le M}j_\ell.
\]
For every $j\ge j_0$ and $i\notin J$, \eqref{fullsquare:local-exclusion} gives $u_{j,i}|_K=0$. Local uniform convergence gives $u_i|_K=0$ as well. Each $i\in J$ has a fixed point $z_i\in V$ with $u_i(z_i)>0$, proving the remaining assertion.
\end{proof}

\begin{lemma}[Uniform cell mass and diameter]\label{fullsquare:diameter}
There exist $\eta_*,m_*,D>0$, depending only on $\Lambda$, such that every reflected minimizing family with $R\ge R_{\rm cmp}(\Lambda)$ satisfies, for every label $i$ and every $x\in D_i$,
\begin{equation}\label{fullsquare:cell-mass}
 \sup_{B_1(x)}u_i\ge\eta_*,
 \qquad \int_{B_2(x)}u_i^2\ge m_*,
 \qquad \operatorname{diam}D_i\le D.
\end{equation}
The mass and amplitude bounds concern the same label $i$, normalized on its complete cell.
\end{lemma}
\begin{proof}
If the first bound failed for this class of families, translations and relabeling would give
\[
 0\in D_{j,i_j},\qquad
 \sup_{B_1}u_{j,i_j}\longrightarrow0,\qquad
 U_j\longrightarrow U
 \quad\text{as in Lemma~\ref{fullsquare:compactness}}.
\]
The limit is nonzero. Apply Lemma~\ref{fullsquare:clearing} with $K=\{0\}$ and $V=B_1$. Eventually $i_j\in J$ for one finite set $J$. For each $i\in J$, choose $y_i\in B_1$ with $u_i(y_i)>0$ and put $\gamma=\min_{i\in J}u_i(y_i)>0$. Local uniform convergence at these finitely many points gives
\[
 \sup_{B_1}u_{j,i_j}\ge u_{j,i_j}(y_{i_j})
 \ge\gamma/2
\]
for large $j$, a contradiction. This proves the uniform constant $\eta_*$.

Let $L$ be the common Lipschitz bound. Choose $y\in B_1(x)$ with $u_i(y)\ge\eta_*/2$, and set
\[
 r_*=\min\{1/2,\eta_*/(4L)\},\qquad
 m_*=\pi r_*^2\eta_*^2/16.
\]
Then $B_{r_*}(y)\subset B_2(x)$ and $u_i\ge\eta_*/4$ on that ball, proving the mass bound.

Join two points of $D_i$ by a path inside the cell. Along the path choose points $x_1,\ldots,x_M$ whose distances from its initial point increase in increments of $5$. Their radius-$2$ balls are disjoint, so complete-cell normalization gives
\[
 1=\int_{\mathbb R^2}u_i^2
 \ge\sum_{k=1}^M\int_{B_2(x_k)}u_i^2
 \ge M m_*.
\]
Thus no such path can join points at distance greater than $5(1+m_*^{-1})$. This value may be taken as $D$.
\end{proof}

\begin{corollary}[Complete mass and the limit network]\label{fullsquare:complete-limits}
Every nonzero label in a limit of Lemma~\ref{fullsquare:compactness} has mass one and compact support of diameter at most $D$. Its approximating labels converge globally in $H^1(\mathbb R^2)$. Each connected component of its canonical positivity set has first eigenvalue $\mu_i$. The limit has a locally finite regular interface network, with pressure matching on regular two-phase arcs and $g\to0$ at degenerate nodes.
\end{corollary}
\begin{proof}
A nonzero limit label has positive mass in a fixed ball. Lemma~\ref{fullsquare:diameter} captures the complete supports of its approximating labels in a larger fixed ball. Local strong $H^1$ convergence therefore becomes global strong $H^1$ convergence, and
\[
 \|u_i\|_2^2=\lim_{j\to\infty}\|u_{j,i}\|_2^2=1,
 \qquad
 \int|\nabla u_i|^2
 =\lim_{j\to\infty}\lambda_{j,i}=\mu_i.
\]
Since $H_0^1(\mathbb R^2)=H^1(\mathbb R^2)$, Lemma~\ref{fullsquare:positivity-trace} gives $H_0^1$ membership on each canonical positivity component $C$. The positive-phase equation and the ground-state identity, for $\phi\in C_c^\infty(C)$,
\[
 \int_C|\nabla\phi|^2-\mu_i\int_C\phi^2
 =\int_C u_i^2|\nabla(\phi/u_i)|^2\ge0,
\]
give $\lambda_1(C)\ge\mu_i$. Testing the equation with the nonzero restriction $u_i\mathbf1_C$ gives the reverse inequality.

Local finiteness now permits application of the finite-phase weak-reflection theorem to the limit. As before, use a connected relatively compact exhaustion containing a fixed nonzero point. The coupling equations, nonnegative defect measures, local Lipschitz continuity, and radial stress identity verify its assumptions. We obtain a locally finite regular network. Stress conservation gives pressure matching on nondegenerate two-phase arcs, while at degenerate nodes $g=|\mathcal Q|\to0$. No connectivity of limit labels or strict inequality $\mu_i<\Lambda$ has been assumed at this stage.
\end{proof}

\section{A positive boundary correction for the square and a smooth example}\label{fullsquare:sec}

Write $Q_R=(0,R)^2$, $B_Q(R)=G(Q_R)+c_*R^2$, and more generally $B_U(R)=B(RU)$. Recall that $\Lambda=2c_*$. Section~\ref{upper:sec:squareupper} gives
\begin{equation}\label{fullsquare:elementary}
 0\le B_Q(R)\le CR\quad(R\ge1),\qquad
 B_Q(\sqrt m)\le m\delta_m(Q).
\end{equation}
We prove a linear lower bound for $B_Q(R)$ by excluding an entire straight interface in a limit satisfying finite replacement comparisons.

\begin{theorem}[Positive square correction]\label{fullsquare:main}
There exist $\sigma>0$ and $R_0<\infty$ such that
\begin{equation}\label{fullsquare:grand-positive}
 B_Q(R)\ge\sigma R\qquad(R\ge R_0).
\end{equation}
Consequently, there exist $\kappa_Q,C_Q>0$ and $m_0\in\mathbb N$ such that
\begin{equation}\label{fullsquare:positive}
 \kappa_Qm^{-1/2}\le\delta_m(Q)\le C_Qm^{-1/2}
 \qquad(m\ge m_0).
\end{equation}
\end{theorem}

\subsection{Finite replacements and absolute defects}
The next proposition passes local minimality to the limit. Its error is the absolute defect of a complete-cell family in its recovery domain.

\begin{proposition}[Finite replacements in a limit]\label{fullsquare:replacement}
Let $U_j\supset B_j$ be bounded open sets containing finite packings $(D_{j,k})_{k\in\mathcal K_j}$ of complete cells from translated reflected configurations with base-square side lengths $R_j\ge R_{\rm cmp}(\Lambda)$. Write $\lambda_{j,k}=\lambda_1(D_{j,k})$ and normalize each eigenfunction on its complete cell. After zero extension and addition of zero coordinates for absent labels, suppose the vectors $\mathbf u_j=(u_{j,k})$ converge as in Lemma~\ref{fullsquare:compactness} to a nonzero limit $U=(u_i)$ described in Corollary~\ref{fullsquare:complete-limits}. Assume
\begin{equation}\label{fullsquare:absolute-defect}
 d_j:=\sum_{k\in\mathcal K_j}(\lambda_{j,k}-\Lambda)-G(U_j)
 \ge0,\qquad d_j\longrightarrow0.
\end{equation}
Let $I$ be any finite set of nonzero limit labels, and let $W_1,\ldots,W_\ell$ be any finite family of bounded connected, pairwise disjoint open sets such that
\[
 u_i|_{W_1\cup\cdots\cup W_\ell}=0\qquad(i\notin I).
\]
The empty insertion family is allowed. Then
\begin{equation}\label{fullsquare:finite-comparison}
 \sum_{i\in I}(\mu_i-\Lambda)
 \le\sum_{\alpha=1}^{\ell}(\lambda_1(W_\alpha)-\Lambda).
\end{equation}
\end{proposition}
\begin{proof}
Fix the finite deletion and insertion families and $\eta>0$, and put $W=\bigcup_\alpha W_\alpha$. In each $W_\alpha$, choose a smooth compactly supported test function, with total Rayleigh error at most $\eta$. Finitely many paths within $W_\alpha$ connect a finite ball covering of its support and give a connected open neighborhood $N_\alpha\Subset W_\alpha$ containing the support. Thus
\[
 \sum_{\alpha=1}^{\ell}\lambda_1(N_\alpha)
 \le\sum_{\alpha=1}^{\ell}\lambda_1(W_\alpha)+\eta.
\]
Choose a smooth $0\le\chi\le1$, zero near all $\overline{N_\alpha}$ and equal to one outside $W$. This cutoff remains fixed while $j\to\infty$. If $\ell=0$, take $\chi=1$.

For the undeleted labels set
\[
 a_j=\sum_{k\in\mathcal K_j\setminus I}\int(1-\chi^2)u_{j,k}^2,
 \qquad
 e_j=\sum_{k\in\mathcal K_j\setminus I}\int|\nabla\chi|^2u_{j,k}^2.
\]
Every undeleted limit coordinate vanishes on $W$, so strong convergence of the full vectors gives
\begin{equation}\label{fullsquare:replacement-error}
 a_j+e_j
 \le (1+\|\nabla\chi\|_\infty^2)
       \|\mathbf u_j-U\|_{L^2(W;\ell^2)}^2
 \longrightarrow0.
\end{equation}
In particular, every retained cutoff has mass at least $1-a_j>0$ for large $j$. Testing the eigenvalue equation with $\chi^2u_{j,k}$ gives, with $\mathcal R$ denoting the Rayleigh quotient,
\[
 \mathcal R(\chi u_{j,k})-\lambda_{j,k}
 =\frac{\int|\nabla\chi|^2u_{j,k}^2}{\int\chi^2u_{j,k}^2},
 \qquad
 \sum_{k\in\mathcal K_j\setminus I}
   \bigl(\mathcal R(\chi u_{j,k})-\lambda_{j,k}\bigr)
 \le\frac{e_j}{1-a_j}.
\]
For each retained label choose one nonzero connected component with quotient at most that of the cutoff, up to an arbitrarily small error. It is disjoint from every $N_\alpha$ and contributes exactly one count term $-\Lambda$.

Corollary~\ref{fullsquare:complete-limits} ensures that every label in $I$ belongs to $\mathcal K_j$ for large $j$, that its complete support lies in a fixed ball, and that $\lambda_{j,i}\to\mu_i$. Eventually all these supports and all $\overline{N_\alpha}$ lie in $B_j\subset U_j$. Delete the labels in $I$, retain the selected cutoff components, and insert the $N_\alpha$. The resulting finite packing has exactly
\[
 n_j^{\rm comp}=\#\mathcal K_j-\#I+\ell
\]
cells. Sending the finitely many component-selection errors to zero, its comparison with $G(U_j)$ is
\[
 G(U_j)\le
 \sum_{k\in\mathcal K_j\setminus I}(\lambda_{j,k}-\Lambda)
 +\sum_{\alpha=1}^{\ell}(\lambda_1(W_\alpha)-\Lambda)
 +\eta+\frac{e_j}{1-a_j}.
\]
Use \eqref{fullsquare:absolute-defect} and cancel the undeleted terms to obtain
\[
 \sum_{i\in I}(\lambda_{j,i}-\Lambda)
 \le\sum_{\alpha=1}^{\ell}(\lambda_1(W_\alpha)-\Lambda)
 +\eta+\frac{e_j}{1-a_j}+d_j.
\]
First let $j\to\infty$, with $\eta$ and $\chi$ fixed, using \eqref{fullsquare:replacement-error} and $d_j\to0$. Then let $\eta\downarrow0$. This proves \eqref{fullsquare:finite-comparison}, including possible labels with $\mu_i=\Lambda$.
\end{proof}

\begin{corollary}[Strict spectral threshold and connected limit cells]\label{fullsquare:limit-cells}
Under the hypotheses of Proposition~\ref{fullsquare:replacement}, every nonzero limit label has $\mu_i<\Lambda$ and a connected canonical positivity set.
\end{corollary}
\begin{proof}
Suppose $\mu_i=\Lambda$. Take a regular adjacency arc of one of its bounded positivity components and enlarge the adjacent component of another label $k$ across a small subarc. The resulting bounded connected open set $W$ occupies only these two phases and their common arc, and
\[
 \lambda_1(W)<\mu_k.
\]
Indeed, equality in domain monotonicity would make the zero extension of the old first eigenfunction a first eigenfunction on $W$, contradicting positivity on its nonempty open enlargement. Deleting the complete labels $i,k$ and inserting only $W$ in \eqref{fullsquare:finite-comparison} gives
\[
 \mu_k-\Lambda
 =(\mu_i-\Lambda)+(\mu_k-\Lambda)
 \le\lambda_1(W)-\Lambda<\mu_k-\Lambda,
\]
a contradiction. Thus $\mu_i<\Lambda$ for every nonzero label.

If a label has two positivity components, each has first eigenvalue $\mu_i$ by Corollary~\ref{fullsquare:complete-limits}. Delete that complete label and insert connected inner approximations to two components. For every $\eta>0$, \eqref{fullsquare:finite-comparison} gives
\[
 \mu_i-\Lambda\le2(\mu_i+\eta-\Lambda).
\]
Taking $0<2\eta<\Lambda-\mu_i$ contradicts this inequality. Hence the positivity set is connected.
\end{proof}

\subsection{Excluding an entire straight interface}
\begin{lemma}[No entire straight interface]\label{fullsquare:no-line}
Suppose a whole-plane configuration has the locally finite regular network and stress conservation described above. Assume that its cells have uniformly bounded diameter and area bounded away from zero, and that each label is a normalized first eigenfunction on its complete connected cell, with eigenvalue $\mu_i<\Lambda$. If \eqref{fullsquare:finite-comparison} holds for every finite deletion and insertion, then its common zero set contains no entire straight line.
\end{lemma}
\begin{proof}
After translation and rotation, suppose that the whole horizontal axis is a common zero set. Connectivity puts each cell entirely in one open half-plane. Choose a maximal connected interval $J$ in the regular part of the axis. Local finiteness of the network ensures that $J$ is nonempty and that the two adjacent labels are fixed on it. Let their complete normalized first eigenfunctions be $u_+,u_-$, with exterior normals $\nu_+=-e_y$ and $\nu_-=e_y$. Stress matching gives
\begin{equation}\label{fullsquare:pressure}
 q(x):=|\partial_{\nu_+}u_+(x,0)|^2
 =|\partial_{\nu_-}u_-(x,0)|^2
 =|\mathcal Q(x,0)|>0,\qquad x\in J.
\end{equation}
The last expression extends the pressure continuously along the axis. Since the complete adjacent cells have diameter at most $D$, the interval has length at most $D$ and finite endpoints. By maximality these endpoints are degenerate nodes, where $|\mathcal Q|\to0$. Thus $q$ cannot be a positive constant on $J$. Choose regular points $x_0,x_1$ with $q(x_0)\ne q(x_1)$ and put $s=x_0-x_1$.

\emph{Tangential sliding.}
Fix $r_0>2D+2|s|+1$, enlarging it so that the required complete cells, the normal-deformation neighborhood, and its translated preimage lie in $B_{r_0}$. For $L>er_0$, let $\ell_L=\log(L/r_0)$ and
\[
 \chi_L(z)=\begin{cases}
 1,&|z|\le r_0,\\
 \log(L/|z|)/\ell_L,&r_0<|z|<L,\\
 0,&|z|\ge L.
 \end{cases}
\]
Apply $\Phi_L(z)=z+s\chi_L(z)e_x$ in the upper half-plane only, leaving the lower half-plane unchanged. For sufficiently large $L$,
\[
 \|D(s\chi_Le_x)\|_\infty\le\frac{|s|}{r_0\ell_L}<\frac12.
\]
The map preserves the vertical coordinate and is strictly increasing on each horizontal line, so it is a bi-Lipschitz self-homeomorphism of the half-plane. Only finitely many complete labels are changed.

Let $X=s\chi_Le_x$ and $\alpha_i=\operatorname*{ess\,sup}_{D_i}|DX|$. Write
\[
 T_i=(|\nabla u_i|^2-\mu_i u_i^2)I-2\nabla u_i\otimes\nabla u_i.
\]
For $\Phi=\operatorname{Id}+X$ and $J_\Phi=\det D\Phi$, the transported test function $u_i\circ\Phi^{-1}$ has quotient
\[
 \mathcal R_i(\Phi)=
 \frac{\displaystyle\int_{D_i}
       |(D\Phi)^{-\mathsf T}\nabla u_i|^2J_\Phi}
      {\displaystyle\int_{D_i}u_i^2J_\Phi}.
\]
Since $\alpha_i<1/2$, expansion of the matrix inverse and determinant gives
\begin{align*}
 \int_{D_i}u_i^2J_\Phi
 &=1+\int_{D_i}u_i^2\operatorname{div}X+r_{i,0},
 & |r_{i,0}|&\le C\alpha_i^2,\\
 \int_{D_i}|(D\Phi)^{-\mathsf T}\nabla u_i|^2J_\Phi
 &=\mu_i+\int_{D_i}
 \bigl(|\nabla u_i|^2\operatorname{div}X
       -2\nabla u_i\cdot DX\nabla u_i\bigr)+r_{i,1},
 & |r_{i,1}|&\le C\mu_i\alpha_i^2.
\end{align*}
The denominator is at least $(1-\alpha_i)^2\ge1/4$. Expanding its reciprocal therefore yields
\begin{equation}\label{fullsquare:transport-quotient}
 \left|\mathcal R_i(\Phi)-\mu_i-\int_{D_i}T_i:DX\right|
 \le C\Lambda\alpha_i^2,
 \qquad \lambda_1(\Phi(D_i))\le\mathcal R_i(\Phi).
\end{equation}
All constants are independent of $i$, the number of labels, and $L$. The term $-\mu_i u_i^2\operatorname{div}X$ comes from the normalization denominator; the bounds use complete-cell normalization and $\int|\nabla u_i|^2=\mu_i\le\Lambda$. Almost-everywhere derivatives suffice.

The sum of the linear terms over all upper cells is zero. For fixed $L$, work in a bounded neighborhood of $\operatorname{supp}X$, remove radius-$\epsilon$ disks around its $N_L$ network nodes, and integrate $\operatorname{div}T_i=0$ in each regular phase. Pressure matching cancels interior-interface fluxes. On the horizontal axis, the eigenfunctions have zero tangential derivative, so the stress flux is normal while $X$ is tangential. On the outer boundary $X=0$. The absolute sum of the node-circle fluxes is bounded by
\[
 C N_L\epsilon\,\|T\|_\infty\|X\|_\infty
 \longrightarrow0\qquad(\epsilon\downarrow0,\ L\text{ fixed}).
\]
The Lipschitz field is justified by convolving $\chi_L$ with compactly supported smooth kernels in the plane and restricting the fields $s\chi_m e_x$ to the upper half-plane. They remain tangential on the axis, and $DX_m\to DX$ locally in $L^1$. Therefore
\[
 \sum_{D_i\subset\{y>0\}}\int_{D_i}T_i:DX=0.
\]

We estimate the quadratic terms using complete cells. Let $d_i=\inf_{z\in D_i}|z|$ and $a_0=b/\Lambda$. First separate those with $d_i<r_0$, some of which cross the boundary of the rigid-translation core. Diameter and area bounds give
\[
 \#\{i:d_i<r_0\}\le\frac{\pi(r_0+D)^2}{a_0},\qquad
 \alpha_i\le\frac{|s|}{r_0\ell_L},\qquad
 \sum_{d_i<r_0}\alpha_i^2\le C\ell_L^{-2}.
\]
Cells contained entirely in $B_{r_0}$ are translated rigidly and have unchanged eigenvalues; the displayed estimate also controls the cells crossing its boundary. Group the other deformed cells into dyadic shells
\[
 2^{k-1}r_0\le d_i<2^kr_0,\qquad k\ge1.
\]
There are $O((2^kr_0)^2/a_0)$ cells in a shell and $\alpha_i\le C|s|/(2^kr_0\ell_L)$. Only $d_i<L$ can contribute, giving $O(\ell_L)$ shells. Hence
\[
 \sum_i\alpha_i^2\le C\ell_L^{-2}+C\ell_L^{-1}\longrightarrow0.
\]
Together with cancellation of the linear terms, the total spectral cost of sliding is at most $C/\ell_L$. The constant depends on the fixed $r_0,D,s,\Lambda,a_0$ and is independent of $L$.

\emph{Normal deformation.}
In a node-free neighborhood of $x_0$, the upper pressure after sliding is $q(x-s)$ and the lower pressure is $q(x)$. The relevant complete upper cell lies in the rigid core and the lower cell has not moved, so their normalized first eigenfunctions are independent of $L$. Choose a smooth common vector field $Y$, compactly supported in this neighborhood, with $Y(x,0)=f(x)e_y$. Apply its small flow $\Psi_t$ to the whole family; only the specified upper and lower cells meet the support. Let $E_{\rm trial}(t)$ be the sum of the normalized Rayleigh quotients of the two transported complete eigenfunctions. Pullback and differentiation of the normalization denominators, followed by local integration by parts, give
\[
 E_{\rm trial}'(0)=\int[q(x-s)-q(x)]f(x)\,dx.
\]
Continuity and $q(x_0-s)\ne q(x_0)$ permit a choice of $f$ with
\[
 \gamma:=-\int[q(x-s)-q(x)]f(x)\,dx>0.
\]
Differentiability at zero gives $t_0>0$, independent of $L$, such that
\[
 E_{\rm trial}(t)-E_{\rm trial}(0)\le-\gamma t/2
 \qquad(0<t<t_0).
\]
Fix one such $t$ and set $\delta=\gamma t/2>0$ before increasing $L$. The deformed first eigenvalues are bounded above by these trial quotients, so the two-cell spectral sum decreases by at least $\delta$. The derivative is computed from the normalized trial quotients and uses integration by parts only near the regular straight arc.

Let $I_+$ contain all complete upper labels meeting the sliding support, and let $I=I_+\cup\{i_-\}$, where $i_-$ is the lower label changed by the local flow. Local finiteness makes $I$ finite. Define
\[
 \widehat D_i=\begin{cases}
 \Phi_L(D_i),&D_i\subset\{y>0\},\\
 D_i,&D_i\subset\{y<0\},
 \end{cases}
 \qquad W_i=\Psi_t(\widehat D_i)\quad(i\in I).
\]
The first map preserves the two half-planes and disjointness; the second is one common plane homeomorphism. Thus the $W_i$ are bounded, nonempty, connected, pairwise disjoint open sets. Every undeleted label is fixed by both maps and vanishes on their union. Delete $I$ and insert $(W_i)_{i\in I}$. The cell count is unchanged, so all $-\Lambda$ terms cancel, and
\[
 \sum_{i\in I}(\lambda_1(W_i)-\Lambda)
 \le\sum_{i\in I}(\mu_i-\Lambda)+C/\ell_L-\delta.
\]
Choosing $L$ with $C/\ell_L<\delta/2$ contradicts \eqref{fullsquare:finite-comparison}.
\end{proof}

\subsection{Seam defects and the square lower bound}
\begin{proof}[Proof of Theorem~\ref{fullsquare:main}]
Reflect an exact minimizing family in a square of side $R\ge R_{\rm cmp}(\Lambda)$ across two adjacent edges to obtain four independent copies in a square of side $2R$. Its auxiliary value is
\[
 4G(RQ)=-4c_*R^2+4B_Q(R).
\]
The preceding lemmas give a uniform diameter bound $D$ and area lower bound $a_0>0$. For a point $y$ on the interior horizontal seam and $r>0$ with $B_{r+D}(y)\Subset Q_{2R}$, put
\[
 I_y(r)=\{i:D_i\cap B_r(y)\ne\varnothing\},\qquad
 U_y(r)=\operatorname{int}\!\bigcup_{i\in I_y(r)}\overline{D_i}.
\]
Then
\begin{equation}\label{fullsquare:recovery}
 \#I_y(r)\le\frac{\pi(r+D)^2}{a_0},\qquad
 B_r(y)\subset U_y(r)\subset B_{r+D}(y).
\end{equation}
For the inner inclusion, the finite closed union cannot miss a point of $B_r(y)$ without missing a nonempty open subset, contrary to almost-everywhere coverage. To verify connectivity, let $E=\bigcup_{i\in I_y(r)}D_i$. Finiteness gives $U_y(r)=\operatorname{int}\overline E$, with $E$ dense in $U_y(r)$. Every selected connected cell belongs to $U_y(r)$ and meets $B_r(y)$, so all lie in the component containing that ball. Any other component would be nonempty and open, hence meet $E$, contradicting connectivity of a selected cell. Thus $U_y(r)$ is connected. An unselected open cell cannot meet the closure of a selected cell: a small ball in the former would then meet the latter, contradicting disjointness. Hence all unselected cells are disjoint from $U_y(r)$.

Define the nonnegative local defect
\[
 \Delta_y(r)=\sum_{i\in I_y(r)}(\lambda_i-\Lambda)-G(U_y(r)).
\]
Only complete cells are removed, with each count term included once. Set $T=r+D+1$ and assume $R\ge16T$. Choose $N=\lfloor R/(8T)\rfloor\ge R/(16T)$ centers, separated by $4T$, on the central half of the interior seam. Their recovery domains and complete label families are disjoint and lie in $Q_{2R}$. Replace each family by a finite packing approaching $G(U_y(r))$. All replacements can be made simultaneously. The large-square bound $G(Q_{2R})\ge-4c_*R^2$, followed by passage of the finitely many errors to zero, gives
\begin{equation}\label{fullsquare:budget}
 \sum_y\Delta_y(r)\le4B_Q(R),\qquad
 \min_y\Delta_y(r)\le64(r+D+1)\frac{B_Q(R)}R.
\end{equation}

Suppose \eqref{fullsquare:grand-positive} fails. Since $B_Q(R)\ge0$, choose $R_j\to\infty$ such that
\[
 R_j\ge\max\{R_{\rm cmp}(\Lambda),16(j+D+1),j^2\},\qquad
 \frac{B_Q(R_j)}{R_j}\le\frac1{64(j+D+1)j^2}.
\]
Use $r=j$ and translate a minimizing center to zero. The recovery domain satisfies $B_j\subset U_j\subset B_{j+D}$, and its absolute defect relative to arbitrary finite connected open packings is at most $j^{-2}$. Lemma~\ref{fullsquare:compactness} and Corollary~\ref{fullsquare:complete-limits} give a nonzero whole-plane limit whose nonzero labels have mass one and support diameter at most $D$. Local uniform convergence preserves the entire horizontal axis as a common zero set.

On every fixed ball the selected family equals the full reflected family for large $j$, since it includes all cells meeting $B_j$. Thus Proposition~\ref{fullsquare:replacement} applies and gives all finite replacement comparisons. Corollary~\ref{fullsquare:limit-cells} excludes threshold eigenvalues and ensures connected limit positivity sets. Lemma~\ref{fullsquare:no-line} now contradicts the retained straight interface. This proves \eqref{fullsquare:grand-positive}.

For each prescribed integer $m$, take $R=\sqrt m$. A dilated original $m$-partition is admissible for $G(RQ)$, so
\[
 B_Q(\sqrt m)\le\frac{S_m(Q)}m-c_*m=m\delta_m(Q).
\]
Therefore $\delta_m(Q)\ge\sigma m^{-1/2}$ for all sufficiently large integers. Proposition~\ref{upper:prop:squareupper} supplies the upper bound.
\end{proof}

\subsection{Comparison with a fixed complement}
To preserve a positive linear correction after smoothing, we estimate the cutoff cost in the removed corner caps. The estimate uses their total boundary length and permits cusps.

\begin{lemma}[Boundary-length comparison for the complement]\label{fullsmooth:transfer}
Let $P$ be bounded and open, and let $\Omega,K\subset P$ be disjoint open sets with $|P\setminus(\Omega\cup K)|=0$ and $\Omega$ connected. Suppose $\partial K$ is covered by $q$ rectifiable curves of total parameter length $L$. Then
\begin{equation}\label{fullsmooth:budget}
 B_\Omega(R)\ge B_P(R)-4\pi(c_*+1)(RL+q).
\end{equation}
\end{lemma}
\begin{proof}
Combining finite packings and using exact area additivity gives $G(RP)\le G(R\Omega)+G(RK)$, hence $B_\Omega(R)\ge B_P(R)-B_K(R)$. In Lemma~\ref{upper:lem:translation}, take $\chi(x)=\min\{1,\operatorname{dist}(x,(RK)^c)\}$. Its positive-level truncations have compact support in $RK$ and converge strongly in $H^1$, proving membership in $H_0^1(RK)$. We obtain
\[
 B_K(R)\le(c_*+1)\bigl|\{x\in RK:
 \operatorname{dist}(x,\partial(RK))<1\}\bigr|.
\]
Sampling a curve of length $\ell$ at arc-length intervals of at most $1$ covers its unit neighborhood by at most $\ell+1$ disks of radius $2$. Summing over the $q$ dilated curves bounds the inner layer by $4\pi(RL+q)$ and proves the lemma. No connectivity or Lipschitz assumption on $K$ is needed.
\end{proof}

\subsection{A fixed smooth example}
\begin{corollary}[A smooth domain attaining the half-order rate]\label{fullsmooth:sharpness}
There is a fixed bounded connected $C^\infty$ domain $\Omega_\varepsilon$ and constants $a_\varepsilon,b_\varepsilon>0$ such that, for every sufficiently large integer $m$,
\[
 a_\varepsilon m^{-1/2}\le\delta_m(\Omega_\varepsilon)
 \le b_\varepsilon m^{-1/2}.
\]
\end{corollary}
\begin{proof}
Let $b(t)=\exp[-1/(t(1-t))]$ for $0<t<1$, with zero endpoint values, and define
\[
 \phi(t)=\frac\pi2\frac{\int_0^t b(s)\,ds}{\int_0^1b(s)\,ds},
 \qquad
 d=\int_0^1\cos\phi(t)\,dt=\int_0^1\sin\phi(t)\,dt\ge\frac12.
\]
The equality follows from symmetry about $1/2$, and the lower bound from $\sin\phi+\cos\phi\ge1$. The corner arc
\[
 \gamma_\varepsilon(t)=(\varepsilon,0)+
 \frac\varepsilon d\int_0^t(-\cos\phi(s),\sin\phi(s))\,ds
\]
joins $(\varepsilon,0)$ to $(0,\varepsilon)$, stays in the corner square, and has one decreasing and one increasing coordinate. Its length is $\varepsilon/d\le2\varepsilon$. Its speed is the positive constant $\varepsilon/d$, and every positive-order endpoint derivative of its tangent angle vanishes. Thus $\gamma_\varepsilon^{(k)}(0)=\gamma_\varepsilon^{(k)}(1)=0$ for $k\ge2$, and joining it to the straight edges in arc-length parametrization gives a regular $C^\infty$ curve.

For $\varepsilon<1/4$, rotate copies into the four disjoint corner squares of $Q$. Together with the remaining straight edges they form a simple closed $C^\infty$ curve bounding a connected domain $\Omega_\varepsilon$. The open complement $K_\varepsilon=Q\setminus\overline{\Omega_\varepsilon}$ has boundary covered by twelve rectifiable pieces of total length at most $16\varepsilon$. The two domain areas sum to $1$.

By \eqref{fullsmooth:budget}, for all sufficiently large real $R$,
\[
 B_{\Omega_\varepsilon}(R)
 \ge[\sigma-64\pi(c_*+1)\varepsilon]R-48\pi(c_*+1).
\]
Fix once and for all
\[
 0<\varepsilon<\min\{1/4,\sigma/[128\pi(c_*+1)]\}.
\]
The domain now remains unchanged as $R$ or $m$ varies, and $B_{\Omega_\varepsilon}(R)\ge\sigma R/4$ for sufficiently large $R$. Let $A=|\Omega_\varepsilon|$ and, for a prescribed large integer $m$, set $R=\sqrt{m/A}$. The original partition competitors give
\[
 m\delta_m(\Omega_\varepsilon)\ge B_{\Omega_\varepsilon}(R)
 \ge\frac{\sigma}{4\sqrt A}\sqrt m.
\]
The upper bound follows from Corollary~\ref{cor:regularity-rates}.
\end{proof}

\begin{proof}[Proof of Theorem~\ref{thm:halfsharp}]
Combine Theorem~\ref{fullsquare:main} and Corollary~\ref{fullsmooth:sharpness}, taking the fixed rounded square as $\Omega_\infty$.
\end{proof}

\section{Optimal rates on H\"older graph domains}\label{holder:sharp-section}

Fix $0<\alpha<1$. We construct one domain $\Omega_\alpha$ and integers $m_j\to\infty$ such that $m_j^{\alpha/2}\delta_{m_j}(\Omega_\alpha)$ stays bounded below by a positive constant. Together with Corollary~\ref{cor:regularity-rates}, this proves Theorem~\ref{thm:holder-sharp}.

The construction is organized around the estimate
\begin{equation}\label{holder:target}
 \delta_{m_j}(\Omega)\ge c h_j-C K_jw_j-C D_jw_j^2/h_j,
 \qquad h_j=\kappa w_j^\alpha,\quad m_j\asymp w_j^{-2},
\end{equation}
proved below in \eqref{holder:delta-ledger}. The core area produces the positive term; the losses come from boundary lengths and the eigenfunction cutoff. Their coefficients $K_j,D_j$ are determined by the prefix. Dirichlet zero traces at adjacent peak heights control the mass on the entire cutoff support, which is needed to estimate the Rayleigh quotients after cutting.

\subsection{Low-energy density and the cost of a cutoff}
Recall the auxiliary packing energy $G$. The perimeter estimate \eqref{upper:eq:perimetergrand} from Section~\ref{upper:sec:squareupper} gives, for any bounded Lipschitz Jordan domain $U$,
\begin{equation}\label{holder:grand-input}
 |G(RU)+c_*|U|R^2|\le C_0(P_U R+1),\qquad R\ge1,
\end{equation}
where $C_0=4\pi(c_*+1)$. This uses only the earlier auxiliary-energy and perimeter bounds, independently of the positive square lower bound in Section~\ref{fullsquare:sec}. We apply this two-sided estimate to the prefix domains and the cores. Every dilated original $m$-partition is also a packing. If $A=|\Omega|$ and $R^2=m/A$, then
\begin{equation}\label{holder:back-to-m}
 \delta_m(\Omega)\ge\frac{G(R\Omega)+c_*AR^2}{AR^2}.
\end{equation}
Neither estimate requires attainment of an infimum.

\begin{lemma}[Low-energy density on arbitrary open cells]\label{holder:heat-density}
Let $D\subset\mathbb R^2$ be bounded and open, and let $u\in H_0^1(D)$ be a normalized first eigenfunction with $\mu=\lambda_1(D)\le\Lambda$. Then
\[
 \|u\|_\infty^2\le e\Lambda/(4\pi)=:M_\Lambda.
\]
For finitely many pairwise disjoint such cells, their zero-extended eigenfunctions satisfy $\sum_i u_i^2\le M_\Lambda$ almost everywhere.
\end{lemma}
\begin{proof}
For arbitrary open $D$, the Dirichlet heat semigroup is dominated by the free heat semigroup; see \cite[Section~1]{Vogt2015} and the exhaustion proof in Appendix~\ref{app:heat}. With the planar heat kernel $k_t(x)=(4\pi t)^{-1}e^{-|x|^2/(4t)}$, \eqref{holder:gaussian-domination} and $T_D(t)u=e^{-\mu t}u$ give
\[
 |u(x)|\le e^{\mu t}\|k_t\|_2\|u\|_2
 =e^{\mu t}(8\pi t)^{-1/2}.
\]
Take $t=(2\Lambda)^{-1}$ and square. At almost every point at most one zero-extended eigenfunction is nonzero, so the same bound holds for their sum of squares, independently of the number of cells.
\end{proof}

\begin{lemma}[A cutoff with control on its entire support]\label{holder:cutoff}
Set $\Lambda=2c_*$ and $M_\Lambda=e\Lambda/(4\pi)$. Let $U$ be bounded and open, let $C\subset\mathbb R^2$ be compact, and let $\chi\in W^{1,\infty}(\mathbb R^2)$ be compactly supported, with $0\le\chi\le1$ and $\chi=1$ on a neighborhood of $C$. Suppose every normalized first eigenfunction of an open cell $D\subset U$ with $\lambda_1(D)\le\Lambda$ satisfies
\begin{equation}\label{holder:full-support}
 \int_{\operatorname{supp}\chi}u^2\le\frac12,
\end{equation}
where $u$ is extended by zero. Then
\[
 G(U)\ge G(U\setminus C)-2M_\Lambda\int_{\mathbb R^2}|\nabla\chi|^2.
\]
\end{lemma}
\begin{proof}
From any finite packing in $U$, delete cells with $\mu_i>2c_*$, whose contributions to $G$ are nonnegative. For each retained cell, put $v_i=(1-\chi)u_i$ and $a_i=\int v_i^2\ge1/2$. Testing with $(1-\chi)^2u_i$ gives
\[
 \frac{\int|\nabla v_i|^2}{a_i}
 =\mu_i+\frac{\int|\nabla\chi|^2u_i^2}{a_i}
 \le\mu_i+2\int|\nabla\chi|^2u_i^2.
\]
Approximate $u_i$ strongly in $H^1$ by $u_{i,l}\in C_c^\infty(D_i)$. Since $\chi=1$ near $C$, $(1-\chi)u_{i,l}$ has compact support in $D_i\setminus C$. Smoothing and passing to the limit show that $v_i\in H_0^1(D_i\setminus C)$.

As in Lemma~\ref{upper:lem:completion}, a Rayleigh quotient is the mass-weighted average of its nonzero component quotients. Select one such component per retained label, with quotient at most the total quotient up to an arbitrarily small error. This preserves the number of retained cells and their disjointness. Lemma~\ref{holder:heat-density} gives
\[
 \sum_i\left(\mathcal R(v_i)-\mu_i\right)
 \le2\int|\nabla\chi|^2\sum_i u_i^2
 \le2M_\Lambda\int|\nabla\chi|^2.
\]
Taking finite nearly minimizing families and sending the finitely many selection and minimization errors to zero proves the claim.
\end{proof}

\subsection{Sparse peak layers on one interval}
Let
\[
 I=[a,b]=[1/4,3/4],\qquad L=b-a=1/2,\qquad\theta=1/256.
\]
Fix a sufficiently small $\kappa>0$. Choose increasing integers $n_j$ and set
\[
 w_j=L2^{-n_j},\qquad h_j=\kappa w_j^\alpha,\qquad
 N_j=L/w_j,\qquad y_{jk}=a+(k+1/2)w_j.
\]
Define triangular layers, extended by zero outside $I$, by
\begin{equation}\label{holder:peaks}
 p_j(y)=h_j\sum_{k=0}^{N_j-1}
 \left(1-\frac{|y-y_{jk}|}{\theta w_j}\right)_+.
\end{equation}
Require $w_{j+1}\le w_j/2$, $h_{j+1}\le h_j/2$, and $2h_1<1/4$. Then $f=\sum_jp_j$ converges uniformly and $0\le f<1/4$. The final domain is
\begin{equation}\label{holder:domain}
 \Omega=\{(x,y):0<y<1,\ f(y)<x<1\},\qquad A=|\Omega|,
\end{equation}
and $3/4<A<1$.

Each zero-extended layer satisfies
\[
 |p_j(y)-p_j(z)|\le
 \min\{\kappa w_j^\alpha,\kappa\theta^{-1}w_j^{\alpha-1}|y-z|\}.
\]
For $s=|y-z|>0$, let $J$ be the largest index with $w_J\ge s$, allowing the first sum below to be empty if none exists. Since $w_j\ge2^{J-j}w_J$ for $j\le J$ and $\alpha-1<0$,
\[
 \sum_{j\le J}w_j^{\alpha-1}s
 \le\frac{s w_J^{\alpha-1}}{1-2^{\alpha-1}}
 \le\frac{s^\alpha}{1-2^{\alpha-1}}.
\]
The remaining wavelengths start below $s$ and decay at least geometrically, so
\[
 \sum_{j>J}w_j^\alpha
 \le\frac{w_{J+1}^\alpha}{1-2^{-\alpha}}
 \le\frac{s^\alpha}{1-2^{-\alpha}}.
\]
If $J$ does not exist, apply this second estimate from $j=1$. Thus
\begin{equation}\label{holder:seminorm}
 |f(y)-f(z)|\le\kappa\left[
 \frac{\theta^{-1}}{1-2^{\alpha-1}}+
 \frac1{1-2^{-\alpha}}\right]s^\alpha.
\end{equation}
The estimate holds across the endpoints of $I$, because the zero-extended layers retain the stated global slope bounds.

Continuity and $0\le f<1/4$ imply
\begin{equation}\label{holder:closure}
 \overline\Omega=
 \{(x,y):0\le y\le1,\ f(y)\le x\le1\}.
\end{equation}
The map
\[
 (s,y)\longmapsto(f(y)+(1-f(y))s,y)
\]
is a homeomorphism from the closed unit square onto \eqref{holder:closure}, carrying its interior onto $\Omega$. Thus $\Omega$ is a Jordan domain. The closure formula also gives $\operatorname{int}\overline\Omega=\Omega$.

The left graph is $C^{0,\alpha}$ by \eqref{holder:seminorm}, and $f$ vanishes near $y=0,1$. The remaining sides are straight, and each corner has the rotated representation $t>|s|$. Hence $\partial\Omega$ admits a finite one-sided $C^{0,\alpha}$ covering as in Definition~\ref{def:graph-windows}, and $|\partial\Omega|=0$.

Write $g_j=\sum_{l<j}p_l$ and $M_j=\operatorname{Lip}(g_j)$. Let $U_{j-1}$ be the prefix domain obtained by replacing $f$ with $g_j$, and $A_{j-1}=|U_{j-1}|$. Each prefix is a polygonal Lipschitz Jordan domain. A layer has area $\theta Lh_j$, so linearity of integration gives
\begin{equation}\label{holder:tail-area}
 A_{j-1}-A=\theta L\sum_{l\ge j}h_l\le2\theta Lh_j.
\end{equation}
\subsection{Zero traces on peak lines and comparison with the core}
Define $I'=[a+L/4,b-L/4]$ and the closed core
\[
 C_j=\{(x,y):y\in I',\ h_j/3\le x-g_j(y)\le h_j/2\}.
\]
Since $g_j+h_j/2<1$ and $I'\Subset(0,1)$, we have $C_j\Subset U_{j-1}$. The core need not lie in the final domain $\Omega$. Its interior is a Lipschitz Jordan domain, with
\begin{equation}\label{holder:core-area}
 |C_j|=Lh_j/12,\qquad
 P_{C_j}\le L\sqrt{1+M_j^2}+h_j/3.
\end{equation}
Fix $0\le\zeta\le1$ in $C_c^\infty((a+L/8,b-L/8))$, equal to one near $I'$, and $0\le\varphi\le1$ in $C_c^\infty((1/8,3/4))$, equal to one near $[1/3,1/2]$. Set
\[
 \eta_j(x,y)=\zeta(y)\varphi((x-g_j(y))/h_j).
\]
The chain rule and integration in $x$ give
\begin{equation}\label{holder:gradient}
 \int|\nabla\eta_j|^2
 \le C\{L(1+M_j^2)/h_j+h_j/L\}\le D'_j/h_j.
\end{equation}
Here $C$ depends only on the fixed cutoffs. For $h_j\le1$, the constant $D'_j$ can be fixed from the prefix, $L$, and the cutoffs before choosing $w_j$. Planar dilation leaves this Dirichlet integral unchanged.

Require $w_j<L/16$ and $M_jw_j\le h_j/8$. For $k=0,\ldots,N_j-2$, between adjacent peak centers put
\[
 X_k=\sup_{[y_{jk},y_{j,k+1}]}g_j+3h_j/4,
 \qquad B_k=(0,X_k)\times(y_{jk},y_{j,k+1}).
\]
At both endpoint heights $y_*$,
\[
 X_k\le g_j(y_*)+7h_j/8<g_j(y_*)+h_j\le f(y_*).
\]
In particular, both horizontal ends lie strictly outside the closure of the final domain:
\begin{equation}\label{holder:exterior-segments}
 (0,X_k)\times\{y_{jk},y_{j,k+1}\}
 \subset\mathbb R^2\setminus\overline\Omega.
\end{equation}
To specify the whole cutoff support, set
\[
 \mathcal L_j=\bigcup_{k=0}^{N_j-1}
       \bigl(\mathbb R\times\{y_{jk}\}\bigr),
 \qquad |\mathcal L_j|=0.
\]
The support of $\zeta$ lies between $a+L/8$ and $b-L/8$, hence strictly between the first and last peak centers because $w_j<L/16$. On a strip between consecutive centers, the support condition for $\varphi$ gives $0<x<g_j(y)+3h_j/4\le X_k$. Therefore
\begin{equation}\label{holder:support-cover}
 \operatorname{supp}\eta_j
 \subset\bigcup_{k=0}^{N_j-2}B_k\ \cup\ \mathcal L_j.
\end{equation}
The exceptional lines have zero area; all other points of the cutoff support belong to the displayed rectangles.

First take $u\in C_c^\infty(R\Omega)$ and extend it by zero as $\widetilde u$. For every $\xi\in(0,RX_k)$, the inequality $f(y_*)-\xi/R>0$ and continuity give an exterior neighborhood of $(\xi,Ry_*)$. Thus $\widetilde u(\xi,\cdot)$ vanishes near both endpoint heights. The one-dimensional Dirichlet Poincar\'e inequality and Fubini's theorem yield
\[
 \int_{RB_k}\widetilde u^2
 \le\frac{R^2w_j^2}{\pi^2}\int_{RB_k}|\partial_y\widetilde u|^2.
\]
The rectangles have disjoint vertical interiors. Summation and \eqref{holder:support-cover} give
\begin{equation}\label{holder:zero-trace}
 \int_{R\operatorname{supp}\eta_j}u^2
 \le\sum_k\int_{RB_k}\widetilde u^2
 \le\frac{R^2w_j^2}{\pi^2}\int_{R\Omega}|\nabla u|^2.
\end{equation}
Smooth compactly supported approximation and strong $H^1(\mathbb R^2)$ convergence of zero extensions extend this estimate to $H_0^1(R\Omega)$, including zero-extended cell eigenfunctions.

If
\begin{equation}\label{holder:mass-condition}
 \Lambda R^2w_j^2/\pi^2\le1/2,
\end{equation}
every normalized first eigenfunction on an open cell $D\subset R\Omega$ with $\lambda_1(D)\le\Lambda$ has mass at most $1/2$ on the entire cutoff support. Lemma~\ref{holder:cutoff}, applied with $\chi=\eta_j(\cdot/R)$ and compact set $RC_j$, gives
\begin{equation}\label{holder:cutting-core}
 G(R\Omega)\ge G(R(\Omega\setminus C_j))-D_j/h_j,
 \qquad D_j=2M_\Lambda D'_j.
\end{equation}

Since $\Omega\subset U_{j-1}$, packings in $R(\Omega\setminus C_j)$ and $R\operatorname{int}C_j$ occupy disjoint subsets of $RU_{j-1}$ and can be combined there. Taking finite nearly minimizing families and then their infima gives
\[
 G(RU_{j-1})\le G(R(\Omega\setminus C_j))+G(R\operatorname{int}C_j).
\]
Use the lower side of \eqref{holder:grand-input} on the prefix, its upper side on the core, and \eqref{holder:cutting-core}. For $R\ge1$ satisfying \eqref{holder:mass-condition}, this yields
\begin{equation}\label{holder:ledger}
 G(R\Omega)+c_*AR^2
 \ge c_*\{|C_j|-(A_{j-1}-A)\}R^2-K_jR-D_j/h_j.
\end{equation}
The net area term $|C_j|-(A_{j-1}-A)$ subtracts the full area loss from the prefix to the final domain. One may take
\[
 K_j=C_0(P_{U_{j-1}}+L\sqrt{1+M_j^2}+1/3+2),
\]
using $h_j\le1$ and $R\ge1$ to absorb constant terms. Thus $K_j,D_j$ are fixed before $w_j$ and need not be bounded uniformly in $j$.

Since $2\theta=1/128$,
\[
 |C_j|-(A_{j-1}-A)\ge Lh_j(1/12-1/128)\ge Lh_j/16.
\]
Set $R^2=m/A$ and apply \eqref{holder:back-to-m} to obtain
\begin{equation}\label{holder:delta-ledger}
 \begin{gathered}
 \delta_m(\Omega)\ge
 \frac{c_*Lh_j}{16A}-\frac{K_j}{AR}-\frac{D_j}{AR^2h_j},\\
 R\ge1,\qquad \Lambda R^2w_j^2/\pi^2\le1/2.
 \end{gathered}
\end{equation}

\subsection{Choosing the scales and fixing the domain}
\begin{proof}[Proof of Theorem~\ref{thm:holder-sharp}]
Fix $q=\min\{1,\pi^2/(16\Lambda)\}>0$. At step $j$, fix the prefix $g_j$ and hence $M_j,K_j,D_j$. Then choose $w_j=L2^{-n_j}$ sufficiently small, with $n_j$ an integer, and set
\[
 h_j=\kappa w_j^\alpha,\qquad m_j=\lfloor q/w_j^2\rfloor.
\]
Impose the geometric and integer conditions
\begin{equation}\label{holder:recursion}
 \begin{aligned}
  &2h_1<1/4,\\
  &w_j<L/16,\qquad M_jw_j\le h_j/8,\qquad
       q/w_j^2\ge2 &&(j\ge1),\\
  &w_j\le w_{j-1}/2,\qquad h_j\le h_{j-1}/2,\qquad
       m_j>m_{j-1} &&(j\ge2),
 \end{aligned}
\end{equation}
and the relative error bound
\begin{equation}\label{holder:payment}
 \frac43K_j\sqrt{2/q}\frac{w_j}{h_j}
 +\frac{8D_j}{3q}\frac{w_j^2}{h_j^2}\le c_*L/32.
\end{equation}
For each fixed prefix, all conditions in \eqref{holder:recursion}--\eqref{holder:payment} hold once the integer $n_j$ is sufficiently large: $\alpha<1$ gives
\[
 \frac{w_j}{h_j}=\kappa^{-1}w_j^{1-\alpha}\longrightarrow0,
 \qquad
 \frac{w_j^2}{h_j^2}=\kappa^{-2}w_j^{2-2\alpha}\longrightarrow0.
\]
Later choices preserve the conditions already imposed. Indeed, nonnegativity gives $f\ge g_j+p_j$, so \eqref{holder:exterior-segments} persists, while geometric height decay gives $\sum_{l\ge j}h_l\le2h_j$, preserving \eqref{holder:tail-area}. The recursion therefore determines a final function $f$ and a single fixed domain $\Omega$. Only after fixing that domain do we test the spectral sequence at the integers $m_j$.

The uniform area bounds imply
\[
 \frac{q}{2w_j^2}\le m_j\le\frac q{w_j^2},\qquad
 R_j^{-1}\le\sqrt{2/q}\,w_j,\qquad
 R_j^{-2}\le2w_j^2/q,\qquad R_j=\sqrt{m_j/A}.
\]
Moreover, $\Lambda R_j^2w_j^2/\pi^2\le4q\Lambda/(3\pi^2)\le1/12$. Hence \eqref{holder:delta-ledger} applies, and \eqref{holder:payment} gives
\[
 \delta_{m_j}(\Omega)\ge c_*Lh_j/32,\qquad
 m_j^{\alpha/2}\delta_{m_j}(\Omega)
 \ge\frac{c_*L\kappa}{32}(q/2)^{\alpha/2}>0.
\]
The full boundary has the required H\"older charts by \eqref{holder:seminorm} and the geometric verification above. The upper bound for all integers makes the upper limit finite, completing the proof.
\end{proof}

The recursion requires $\alpha<1$: at $\alpha=1$, $w_j/h_j$ is constant and cannot absorb the changing prefix errors. The square argument supplies that endpoint. For related multiscale graph constructions for the Neumann counting function, see \cite[Theorem~1.10 and \S4.5]{NetrusovSafarov2005}.

\section{Little H\"older domains and common improvement factors}\label{holder:little-section}

\begin{proof}[Proof of Theorem~\ref{thm:little-holder}]
Fix a little H\"older domain and a finite covering as in Definition~\ref{def:little-holder-domain}. Let $\omega$ be the maximum of the moduli of continuity of its graph functions. Then $\omega(t)/t^\alpha\to0$, and \eqref{upper:eq:graphmodulus} gives
\[
 |\delta_m(\Omega)|\le C_\Omega
 \{m^{-1/2}+\omega(K_\Omega m^{-1/2})\}.
\]
After multiplication by $m^{\alpha/2}$, the first term tends to zero since $\alpha<1$, and the second is a fixed multiple of $\omega(t)/t^\alpha$, with $t=K_\Omega m^{-1/2}\to0$. This proves the first assertion.

Now prescribe any positive sequence $\varepsilon_m\to0$. Use the preceding peak layers, cores, and fixed cutoffs, but set
\[
 a_j=2^{-j},\qquad h_j=\kappa a_jw_j^\alpha.
\]
At step $j$, fix the prefix and $a_j>0$ before choosing $w_j$. Impose the full system \eqref{holder:recursion} and the relative error bound \eqref{holder:payment} with these new heights, as well as
\[
 m_j=\lfloor q/w_j^2\rfloor,\qquad \varepsilon_{m_j}\le a_j/j.
\]
To impose the last condition, use $\varepsilon_m\to0$ to choose $T_j\in\mathbb N$ such that
\[
 \varepsilon_m\le a_j/j\qquad\text{for every integer }m\ge T_j,
\]
and then require $\lfloor q/w_j^2\rfloor\ge T_j$. This uses no monotonicity of $\varepsilon$. The two relative errors in \eqref{holder:payment} are now
\[
 O_j(a_j^{-1}w_j^{1-\alpha})\quad\text{and}\quad
 O_j(a_j^{-2}w_j^{2-2\alpha}),
\]
which still tend to zero for the fixed current $a_j$. Halving of both heights and wavelengths, and $M_jw_j\le h_j/8$, can be imposed simultaneously. These choices determine a single domain $\Omega_\varepsilon$, depending on the sequence $\varepsilon$ but fixed before taking the spectral limit.

We verify the little H\"older condition on the final boundary. For fixed $J$, the first $J-1$ layers form a Lipschitz function, whose modulus divided by $t^\alpha$ tends to zero. For each tail layer,
\[
 |p_j(y)-p_j(z)|\le\kappa\Bigl(\sup_{l\ge J}a_l\Bigr)
 \min\{w_j^\alpha,\theta^{-1}w_j^{\alpha-1}|y-z|\}.
\]
The two-part summation used in \eqref{holder:seminorm} gives
\[
 \omega_{{\rm tail},J}(t)\le
 C_\alpha\kappa\sup_{j\ge J}a_j\,t^\alpha
 =C_\alpha\kappa2^{-J}t^\alpha,
\]
where $C_\alpha$ depends only on the fixed $\alpha$ and $\theta$. For each fixed $J$, writing $M_J$ for the Lipschitz constant of the first $J-1$ layers gives
\[
 \frac{\omega_f(t)}{t^\alpha}
 \le M_Jt^{1-\alpha}+C_\alpha\kappa2^{-J},
 \qquad
 \limsup_{t\downarrow0}\frac{\omega_f(t)}{t^\alpha}
 \le C_\alpha\kappa2^{-J}.
\]
First $t\downarrow0$ and then $J\to\infty$ prove $\omega_f(t)=o(t^\alpha)$, including across the endpoints of $I$. The remaining charts are Lipschitz and hence little $C^{0,\alpha}$. The preceding homeomorphism and chart representations verify the Jordan property and the finite one-sided covering; \eqref{holder:closure} gives $\Omega_\varepsilon=\operatorname{int}\overline{\Omega_\varepsilon}$.

For the modified heights, $f\ge g_j+p_j$ and $\sum_{l\ge j}h_l\le2h_j$ preserve \eqref{holder:exterior-segments} and \eqref{holder:tail-area}. The proof of \eqref{holder:delta-ledger} therefore applies, with $K_j,D_j$ determined by the prefix and fixed cutoffs before choosing $w_j$. Set $A=|\Omega_\varepsilon|$ and $R_j=\sqrt{m_j/A}$. The area bounds $3/4<A<1$ give the same estimates for $R_j$ as in Section~\ref{holder:sharp-section}, so \eqref{holder:payment} yields
\[
 \delta_{m_j}(\Omega_\varepsilon)
 \ge\frac{c_*L\kappa}{32}a_jw_j^\alpha
 \ge d_\alpha a_jm_j^{-\alpha/2},\qquad
 d_\alpha=\frac{c_*L\kappa}{32}(q/2)^{\alpha/2}>0.
\]
Consequently,
\[
 \frac{\delta_{m_j}(\Omega_\varepsilon)}
      {\varepsilon_{m_j}m_j^{-\alpha/2}}
 \ge d_\alpha\frac{a_j}{\varepsilon_{m_j}}
 \ge d_\alpha j\longrightarrow+\infty.
\]
The first part of the proof also applies to this fixed domain and gives $|\delta_m(\Omega_\varepsilon)|=o(m^{-\alpha/2})$ on the full integer sequence.
\end{proof}

\subsection{Further questions}
The results determine optimal exponents over the stated domain classes through fixed examples. They leave open the existence of a positive matching half-order lower bound on every fixed Lipschitz or smooth domain, including the disk. A further question is whether $\sqrt m\,\delta_m(\Omega)$ has a limit on suitable domains and how such a coefficient depends on the boundary. Identification of the bulk constant $c_*$ with the conjectured honeycomb value is a separate problem; none of the arguments here requires that identification.

\appendix
\section{Heat-semigroup domination on arbitrary open sets}\label{app:heat}

We record the exhaustion argument used in Lemma~\ref{holder:heat-density}. It applies to the Dirichlet form on $H_0^1(D)$ without any regularity assumption on the cell boundary.

\begin{lemma}\label{app:heat-domination}
Let $D\subset\mathbb R^2$ be bounded and open, and let $T_D(t)$ be the heat semigroup defined by the Dirichlet form on $H_0^1(D)$. For every $v\in L^2(D)$ and $t>0$, with zero extension $\widetilde v$,
\begin{equation}\label{holder:gaussian-domination}
 |T_D(t)v|(x)\le(k_t*|\widetilde v|)(x)
 \quad\text{for almost every }x\in D,\qquad
 k_t(x)=(4\pi t)^{-1}e^{-|x|^2/(4t)}.
\end{equation}
\end{lemma}
\begin{proof}
Take an increasing smooth inner exhaustion $D_l\Subset D$, eventually containing every compact subset of $D$; finite unions of smooth sets are allowed. Let $P_l:L^2(D)\to L^2(D_l)$ be restriction and $J_l:L^2(D_l)\to L^2(D)$ be zero extension. For nonnegative initial data, comparison on smooth domains, positivity, and $L^2$ approximation yield
\[
 0\le J_lT_{D_l}(t)P_lv
 \le k_t*\widetilde{P_lv}\le k_t*\widetilde v.
\]
Here $\widetilde{P_lv}$ is the extension from $D_l$ to the plane.

To pass to the limit, let $A_l,A$ be the nonnegative Dirichlet operators and set
\[
 K_l=(A_l+1)^{-1},\qquad R_l=J_lK_lP_l,\qquad R=(A+1)^{-1}.
\]
The increasing union of the zero-extended spaces $H_0^1(D_l)$ is dense in $H_0^1(D)$ because it contains every smooth compactly supported test function eventually. For $v\in L^2(D)$, $R_lv$ and $Rv$ minimize
\[
 \int_D(|\nabla w|^2+|w|^2)-2\operatorname{Re}\int_Dv\overline w
\]
over $H_0^1(D_l)$ and $H_0^1(D)$, respectively. This is orthogonal projection in the energy inner product, so $R_lv\to Rv$ strongly in $H^1$ and hence in $L^2$.

For fixed $t>0$, define $F_t(s)=\exp[-t(s^{-1}-1)]$ for $s>0$ and $F_t(0)=0$. It is continuous on $[0,1]$ and can be approximated uniformly by polynomials with zero constant term. Since $P_lJ_l=I$, we have $R_l^k=J_lK_l^kP_l$ for $k\ge1$, hence $p(R_l)=J_lp(K_l)P_l$ when $p(0)=0$. Uniform approximation gives
\[
 F_t(R_l)=J_lF_t(K_l)P_l=J_lT_{D_l}(t)P_l.
\]
The bounds $\|R_l\|,\|R\|\le1$ and strong convergence of $R_l$ imply convergence of the polynomials, and then of the continuous functional calculus. Therefore
\[
 J_lT_{D_l}(t)P_lv\longrightarrow T_D(t)v
 \quad\text{strongly in }L^2(D).
\]
An almost-everywhere convergent subsequence preserves comparison; positivity treats general initial data and proves \eqref{holder:gaussian-domination}.
\end{proof}


\begin{thebibliography}{99}
\raggedright

\bibitem{ArmitageGardiner2001}
D.~H. Armitage and S.~J. Gardiner,
\emph{Classical Potential Theory},
Springer Monographs in Mathematics, Springer, London, 2001.

\bibitem{BucurButtazzoHenrot1998}
D. Bucur, G. Buttazzo, and A. Henrot,
\emph{Existence results for some optimal partition problems},
Adv. Math. Sci. Appl. \textbf{8} (1998), no.~2, 571--579.

\bibitem{CaffarelliLin2007}
L.~A. Caffarelli and F.-H. Lin,
\emph{An optimal partition problem for eigenvalues},
J. Sci. Comput. \textbf{31} (2007), 5--18.
\url{https://doi.org/10.1007/s10915-006-9114-8}.

\bibitem{CaffarelliLin2008}
L.~A. Caffarelli and F.-H. Lin,
\emph{Singularly perturbed elliptic systems and multi-valued harmonic functions with free boundaries},
J. Amer. Math. Soc. \textbf{21} (2008), no.~3, 847--862.
\url{https://doi.org/10.1090/S0894-0347-08-00593-6}.

\bibitem{FrankLarsonError2020}
R.~L. Frank and S. Larson,
\emph{On the error in the two-term Weyl formula for the Dirichlet Laplacian},
J. Math. Phys. \textbf{61} (2020), 043504.
\url{https://doi.org/10.1063/1.5145003}.

\bibitem{FrankLarson2020}
R.~L. Frank and S. Larson,
\emph{Two-term spectral asymptotics for the Dirichlet Laplacian in a Lipschitz domain},
J. Reine Angew. Math. \textbf{766} (2020), 195--228.
\url{https://doi.org/10.1515/crelle-2019-0019}.

\bibitem{FrankLarson2025}
R.~L. Frank and S. Larson,
\emph{Riesz means asymptotics for Dirichlet and Neumann Laplacians on Lipschitz domains},
Invent. Math. \textbf{241} (2025), 999--1079.
\url{https://doi.org/10.1007/s00222-025-01352-x}.

\bibitem{GengLin2022}
Z. Geng and F. Lin,
\emph{Large $m$ asymptotics for minimal partitions of the Dirichlet eigenvalue},
Sci. China Math. \textbf{65} (2022), 1--8.
\url{https://doi.org/10.1007/s11425-020-1802-6}.

\bibitem{Henrot2006}
A. Henrot,
\emph{Extremum Problems for Eigenvalues of Elliptic Operators},
Frontiers in Mathematics, Birkh\"auser, Basel, 2006.
\url{https://doi.org/10.1007/3-7643-7706-2}.

\bibitem{NetrusovSafarov2005}
Yu. Netrusov and Yu. Safarov,
\emph{Weyl asymptotic formula for the Laplacian on domains with rough boundaries},
Comm. Math. Phys. \textbf{253} (2005), 481--509.
\url{https://doi.org/10.1007/s00220-004-1158-8}.

\bibitem{Stein1993}
E.~M. Stein,
\emph{Harmonic Analysis: Real-Variable Methods, Orthogonality, and Oscillatory Integrals},
with the assistance of T.~S. Murphy,
Princeton Mathematical Series 43, Princeton University Press, Princeton, NJ, 1993.

\bibitem{TavaresTerracini2012}
H. Tavares and S. Terracini,
\emph{Regularity of the nodal set of segregated critical configurations under a weak reflection law},
Calc. Var. Partial Differential Equations \textbf{45} (2012), 273--317.
\url{https://doi.org/10.1007/s00526-011-0458-z}.

\bibitem{Vogt2015}
H. Vogt,
\emph{$L_1$-estimates for eigenfunctions and heat kernel estimates for semigroups dominated by the free heat semigroup},
J. Evol. Equ. \textbf{15} (2015), no.~4, 879--893.
\url{https://doi.org/10.1007/s00028-015-0285-3}.

\end{thebibliography}
\end{document}